\magnification=1200

\font\twelverm=cmr12 at 14pt
\hsize = 31pc
\vsize = 44pc
\overfullrule = 0pt

\input amssym.def
\input amssym.tex
 \font\newrm =cmr10 at 24pt
\def\bul{\raise .9pt\hbox{\newrm .\kern-.105em } }

 \def\fr{\frak}

 \baselineskip 20pt
 
 \def\h{\hbox{ }}

 \def\u{{\fr u}}
 \def\r{{\fr r}}

 \def\n{{\fr n}}
 \def\a{{\fr a}}
 \def\d{{\fr d}}
 
 \def\ss{{\fr s}}
 
 \def\b{{\fr b}}

 \def\ee{{\fr e}}

 \def\<{\le}
 \def\>{\ge}

 \def\s{{\h\subset\h}}
 
 \def\vs{\vskip }

 \def\mapright#1
  {\smash{\mathop
  {\longrightarrow}
  \limits^{#1}}}

 \def\kk#1{{\kern .4 em} #1}
 \def\vs{\vskip 1pc}

\rm

\font\smallbf=cmbx10 at 9pt 
\font\twelverm=cmr12 at 14pt

\font\authorfont=cmr10 at 11pt
\font\ninerm=cmr9
\rm
\centerline{\twelverm Gelfand--Zeitlin theory from the perspective of classical
 mechanics  II}\vskip 1.2pc 
 \baselineskip=11pt
\vskip8pt
\centerline{\authorfont BERTRAM
KOSTANT\footnote*{\ninerm Research supported in part by NSF grant
DMS-0209473 and in part by the \hfil\break KG\&G Foundation} and
NOLAN WALLACH\footnote{**}
{\ninerm
Research supported in part by NSF grant MTH 0200305}}
\vskip 2pc\baselineskip=11pt 
{\font\ninerm=cmr9 \baselineskip=14pt 
\noindent{\smallbf ABSTRACT.  }\ninerm \rm In this paper, Part II, of a two part paper
we apply the results of [KW], Part I, to establish, with an explicit dual coordinate system, a commutative
analogue of the Gelfand-Kirillov theorem for $M(n)$, the algebra of $n\times n$ complex matrices. The
function field
$F(n)$ of
$M(n)$ has a natural Poisson structure and an exact analogue would be to show that $F(n)$ is isomorphic
to the function field of a $n(n-1)$-dimensional phase space over a Poisson central rational function
field in $n$ variables. Instead we show that this the case for a  Galois extension, $F(n,\ee)$,
of $F(n)$. The techniques uses a maximal Poisson commutative algebra of functions arising from
Gelfand-Zeitlin theory, the algebraic action of a $n(n-1)/2$--dimensional torus on $F(n,\ee)$,
and the structure of a Zariski open subset of $M(n)$ as a $n(n-1)/2$--dimensional torus bundle over a
$n(n+1)/2$--dimensional base space of Hessenberg matrices.

 \vskip 1.5pc \centerline{\bf 0. Part II continuation of Introduction}\vskip 1.5pc 
0.6. We recall some of the notation and
results in Part I, i.e., [KW]. If $k$ is a positive integer then $I_k= \{1,\ldots,k\}$. $M(n)$ is the algebra of all $n\times
n$ complex matrices. If $m\in I_n$, then regard $M(m)\s M(n)$ as the upper left block of all $m\times m$ matrices. If $x\in
M(n)$ then $x_m\in M(m)$ is the upper left principal $m\times m$ minor of $x$. Using a
natural isomorphism of
(the Lie algebra) $M(n)$ with its dual space, $M(n)$ becomes a Poisson manifold so that its affine ring
${\cal O}(M(n))$ is a Poisson algebra. For any
$k\in
\Bbb Z_+$ let
$d(k) = k(k+1)/2$. The subalgebra
$J(n)$ of
${\cal O}(M(n))$, generated by the symmetric polynomial $Gl(m)$-invariants of $M(m)$ for  all $m\in I_n$, is a polynomial algebra
with $d(n)$ generators and more importantly it is a maximal Poisson commutative subalgebra of ${\cal O}(M(n))$. 

In Part I we
showed that the Poisson vector field $\xi_f$ on $M(n)$ corresponding to any $f\in J(n)$ is globally integrable on $M(n)$ and a
choice of generators of $J(n)$ defines an abelian Lie group $A$ of dimension $d(n-1)$ operating on $M(n)$. The orbits of $A$
are explicitly determined in Part I and the orbits are independent of the choice of generators. One particular choice are the
functions, $p_i(x),\,i\in I_{d(n)},\,x\in M(n)$, where, for all $m\in I_n$, $p_{d(m-1)+ k}(x),\,k\in I_m$, are the nontrivial
coefficients of the  characteristic polynomial of $x_m$. 

A suitable measure on
$\Bbb R$ and the Gram-Schmidt process define a sequence, $\phi_k(t)$,\break $k\in \Bbb Z_+$, of
orthogonal polynomials on $\Bbb R$. Let $W_n$ be the span of $\phi_{m-1},\,m\in I_n$, and let
$x\in M(n)$ be the matrix, with respect to this basis, of the operator of multiplication by $t$,
followed by projection on $W_n$. The matrix $x$ is Jacobi and for $m\in I_n$ one recovers the
orthogonal polynomial $\phi_m$ as the characteristic polynomial of $x_m$. In particular the all
important zeros of the orthogonal polynomials $\phi_m$ appear as the eigenvalues of the $x_m$.
One motivation for our work here is to set up Poisson machinery to deal with the eigenvalues of
$x_m$ for any $x\in M(n)$. In the course of setting up this machinery we have obtained a number
of new results. Some of these results have appeared in Part I ([KW]). In the present paper, Part
II, we will be concerned with establishing a refinement of a commutative analogue of the
Gelfand-Kirillov theorem. The refinement refers to exhibiting an explicit coordinate system
satisfying the Poisson commutation relations of phase space. The coordinate system emerges from
the action of an algebraic group and the structure of $M(n)$, obtained in Part I, as sort of a
cotangent bundle over the variety of Hessenberg matrices. 

In more detail let $M_{\Omega}(n)$ be
the Zariski open (dense) subset of $M(n)$ defined as the set of all $x\in M(n)$ such that $x_m$
is regular semisimple in $M(m)$ for all $m\in I_n$ and such that the spectrum of $x_{m-1}$, for
$m>0$, has empty intersection with the spectrum of $x_m$. In Part I $M_{\Omega}(n)$ was shown to
have the following structure: It is a $(\Bbb C^{\times})^{d(n-1)}$ bundle over a
$d(n)$-dimensional base space ($d(n-1) + d(n) = n^2$) base space. The fibers are not only the
level sets for the functions in $J(n)$ but also the fibers are the orbits of $A$ in
$M_{\Omega}(n)$. The base space, denoted by $\b_{e,\Omega(n)}$, is the intersection $\b_e\cap
M_{\Omega}(n)$ where $\b_e$ is the space of Hessenberg matrices. That	is, $x\in \b_e$ if and
only if $x$ is of the form $$x = \left(\matrix{a_{1\,1}&a_{1\,2}&\cdots &a_{1\,n-1}&a_{1\,n}\cr
1&a_{2\,2}&\cdots &a_{2\,n-1}&a_{2\,n}\cr 0&1&\cdots &a_{3\,n-1}&a_{3\,n}\cr \vdots &\vdots
&\ddots &\vdots &\vdots\cr 0&0&\cdots &1&a_{n\,n}\cr}\right)$$ Also $\b_{e,\Omega(n)}$ is
Zariski dense in $\b_e$. Theorem 2.5 in Part I concerning establishing a beautiful property of
$\b_e$ plays a major role here. In part II a ``Lagrangian" property of $\b_e$ plays a key role in
showing the dual coordinates $s_i$, defined below, Poisson commute. Part I has 3 sections. 
A
serious deficiency in $M_{\Omega}(n)$ in dealing with the commutative analogue of the
Gelfand-Kirillov theorem is that one cannot consistently solve the chararteristic polynomials
of $x_m$, for all $m\in I_n$ and all $x\in M_{\Omega}(n)$ to yield algebraic eigenvalue
functions $r_i$ on $M_{\Omega}(n)$. 
\vs 0.7. We begin, in the first section of Part II, labeled
Section 4, to obtain such functions on a covering space $M_{\Omega}(n,\ee)$ of $M_{\Omega}(n)$.
Initially the covering map $$\pi_n: M_{\Omega}(n,\ee)\to M_{\Omega}(n)$$ is only understood to
be analytic. The covering admits, as deck transformations, a group, $\Sigma_n$, isomorphic to
the direct product of the symmetric groups $S_m,\,m\in I_n$ and as analytic manifolds
$$M_{\Omega}(n,\ee)/ \Sigma_n\cong M_{\Omega}(n)$$ However much more structure is needed and
established in \S 4. For one thing $M_{\Omega}(n,\ee)$ is a nonsingular affine variety and
$\pi_n$ is a finite \'etale morphism. For another, if $F(n)$ is the field of rational functions
on $M(n)$ and $F(n,\ee)$ is the field of rational functions on $M_{\Omega}(n,\ee)$, then
$F(n,\ee)$ is a Galois extension of $F(n)$ with $\Sigma_n$ as Galois group. Furthermore the
affine ring ${\cal O}(M_{\Omega}(n,\ee))$ is the integral closure of ${\cal O}(M_{\Omega})$ in
$F(n,\ee)$. Very significant for our purposes, there exists (eigenvalue) functions $r_i\in {\cal
O}(M_{\Omega}(n,\ee)),\,i\in I_{d(n)}$ with the property that for any $m\in I_n$ and any $z\in
M_{\Omega}(n,\ee)$ the numbers $r_{d(m-1) + k}(z),\,k\in I_m$ are the eigenvalues of $x_m$,
where $x=\pi_n(z)$. The Poisson structure on $M_{\Omega}(n)$ lifts to $M_{\Omega}(n,\ee)$ and
one has $[r_i,r_j] = 0$ for all $i,j\in I_{d(n)}$.
\vs 0.8. \S 5 is devoted to the construction of
the dual coordinates $s_j\in {\cal O}(M_{\Omega}(n,\ee)),\,j\in I_{d(n-1)}$. There are two key
points here. (a) It is shown that the Poisson vector fields $\xi_{\r_i}$ on $M_{\Omega}(n,\ee)$
integrate and generate a complex algebraic torus, $A_{\r} \cong (\Bbb C^{\times})^{d(n-1)}$
which operates algebraically on $M_{\Omega}(n,\ee))$ and in fact if $M_{\Omega}(n,\ee,\b)$ is
the $\pi_n$ inverse image of $\b_{\ee,\Omega(n)}$ in $M_{\Omega}(n,\ee)$, then the map
$$A_{\r}\times M_{\Omega}(n,\ee,\b)\to M_{\Omega}(n,\ee),\qquad (b,y) \mapsto \b\cdot y $$ is an
algebraic isomorphism. The natural coordinate system on $A_{\r}$ then carries over to
$M_{\Omega}(n,\ee)$ defining functions $s_j\in {\cal O}(M_{\Omega}(n,\ee)),\,j\in I_{d(n-1)}$,
when they are normalized so that, for all $j$, $s_j$ is the constant 1 on
$M_{\Omega}(n,\ee,\b)$. The second key point, (b), yields the Poisson commutativity $[s_i,s_j] =
0$ from the Lagrangian property of $\b_e$. See Theorem 5.20 and its proof. Combining Theorems
5.14 and 5.23 one has \vs {\bf Theorem 0.16.} {\it The image of the map $$M_{\Omega}(n,\ee)\to
\Bbb C^{n^2},\qquad z \mapsto (r_1(z),\ldots,r_{d(n)}(z),s_1(z),\ldots,s_{d(n-1)}(z)) \eqno
(0.13)$$ is a Zariski open set $Y$ in $\Bbb C^{n^2}$ and (0.13) is an algebraic isomorphism of
$M_{\Omega}(n,\ee)$ with $Y$. Furthermore one has the following Poisson commutation relations:
$$ \eqalign{(1)&\,\,[r_i,r_j]= 0,\,i,j\in I_{d(n)}\cr (2)&\,\,[r_i,s_j] =
\delta_{i\,j}\,s_j,\,i\in I_{d(n)},\,j\in I_{d(n-1)}\cr (3)&\,\,[s_i,s_j] = 0,\,i,j\in
I_{d(n-1)}\cr}\eqno (0.14)$$} \vskip .5pc Noting that $s_i$ vanishes nowhere on
$M_{\Omega}(n,\ee)$ one has $r_{(i)}\in {\cal O}(M_{\Omega}(n,\ee))$ for $i\in I_{d(n-1)}$ where
$r_{(i)} = r_i/s_i$. Replacing $r_i$ by $r_{(i)}$ in (2) one has the more familiar phase space
commutation relation $[r_{(i)},s_j] = \delta_{i\,j}$. For the implication of Theorem 0.16 on the
structure of the field $F(n,\ee)$ see Theorem 5.24. \vs 0.9. Of course given the eigenvalue
functions $r_i$, the dual coordinates $s_j$ are not uniquely determined. In the present paper
they are given by the use of the algebraic group $A_{\r}$ (defined by the $r_i$) and a set of
Hessenberg matrices as a base space. Independently and quite differently the papers [GKL1] and
[GKL2] also deal with establishing a refined commutative analogue of the Gelfand-Kirillov
theorem. A point of similarity is the use of the coordinates $r_i$ and the necessity, thereby,
to go to a a covering. In \S 3 of [GKL1] dual coordinates are given, denoted in that paper by
$Q_{n\,j}$. It seems to be an interesting question to write down equations expressing a relation
between the $Q_{n\,j}$ in [GKL1] and the $s_i$ here. The first three sections of the paper are in
Part I. Part II begins with \S 4.

	\vskip 1.2pc
\centerline{\bf 4. The covering $M_{\Omega}(n,\ee)$ of $M_{\Omega}(n)$ and the eigenvalue functions
$r_i$}
\vskip 1.2pc 4.1. We retain the general notation of Part 1 so that
$n$ is a positive integer and $M(n)$ is the space of all
complex $n\times n$ matrices. As in (2.61),  for $m\in I_n$, (see \S
1.1) let
$\d(m)$ be the space of all diagonal matrices 
in $M(m)$ and let $\ee(m)$ be the 
(connected) Zariski open subset of all regular elements in $\d(m)$. 
That is, if $z \in \d(m)$, then $z\in
\ee(m)$ if and only if the diagonal entries of $z$ are distinct. Consider the direct product
$$\ee =\ee(1)\times \cdots \times \ee(m)\eqno (4.1)$$ so that if $\nu\in \ee$ we can write
$$\nu= (\nu(1),\ldots,\nu(n))\eqno (4.2)$$ where $\nu(m)\in \ee(m)$. In addition we will write
$$ \nu(m) = diag(\nu_{1\,m},\ldots,\nu_{m\,m})\eqno (4.3)$$  where the numbers
$\nu_{i\,m}\in \Bbb C,\,i\in I_m$, are distinct. Taking notation from (2.53) let
$\ee_{\Omega(n)}$ be the Zariski open subset of $\ee$ defined so that if $\nu\in \ee$
then $\nu\in \ee_{\Omega(n)}$  if and only if
$$\nu_{i\,m}\neq \nu_{j\,m+1},\,\,\forall m\in I_{n-1},\,i\in I_m,\,j\in I_{m+1}\eqno (4.4)$$
Of course $\ee_{\Omega(n)} $ is a nonsingular variety where $$dim\,\ee_{\Omega(n)} = d(n)$$ (see
\S 0.1). The symmetric group $S_m$, as the Weyl group of $(M(m),\d(m))$, operates freely on
$\ee(m)$ and the direct product $\Sigma_n = S_1\times \cdots\times S_n$ ( a group of
order $\prod_{m\in I_n}\,m!$) operates freely on
$\ee_{\Omega(n)}$ where if $\sigma = (\sigma_1,\ldots,\sigma_n),\,\,\sigma_m\in S_m$, is in
$\Sigma_n$ and
$\nu\in \ee_{\Omega(n)}$ then, using the notation of (4.2) and (4.3), $$\sigma\cdot\nu = 
(\sigma_1\cdot\nu(1)),\ldots,\sigma_m\cdot \nu(m))\eqno (4.5)$$ and $$\sigma_m\cdot \nu(m) =
diag\,(\nu_{\sigma_m^{-1}(1)\,m},\ldots,\nu_{\sigma_m^{-1}(m)\,m})\eqno (4.6)$$\vskip .5pc 
4.2. Recall the Zariski open set $M_{\Omega}(n)$ of $M(n)$ (see (2.53)). In particular
we recall that the matrices in $M_{\Omega}(n)$ are regular semisimple. Consider the direct
product $\ee_{\Omega(n)}
\times M_{\Omega}(n)$ and let $$M_{\Omega}(n,\ee)= \{(\nu,x)\in \ee_{\Omega(n)}
\times M_{\Omega}(n)\mid \nu(m)\,\,\hbox{is $Gl(m)$-conjugate to $x_m$},\,\,\forall m\in
I_n\}\eqno (4.7)$$ It is clear that $M_{\Omega}(n,\ee)$ is a Zariski closed subset of
$\ee_{\Omega(n)} \times M_{\Omega}(n)$ and the maps $$\pi_n:M_{\Omega}(n,\ee)\to 
M_{\Omega}(n)\,\,\,\,\hbox{where}\,\,\pi_n(\nu,x) = x\eqno (4.8)$$ and $$\kappa_n:
M_{\Omega}(n,\ee)\to \ee_{\Omega(n)}\,\,\,\,\hbox{where}\,\,\kappa_n(\nu,x)=\nu\eqno
(4.9)$$ are surjective (see Theorem 2.5) algebraic morphisms. 

For $m\in I_n$ and $i\in I_m$ let $\rho_{i\,m}$ be the regular function on
$M_{\Omega}(n,\ee)$ defined so that if $z\in M_{\Omega}(n,\ee)$ and $\nu = \kappa_n(z)$ then
$\rho_{i\,m}(z) = \nu_{i\,m}$. \vs {\bf Remark 4.1.} One notes that if $z\in M_{\Omega}(n,\ee)$
and $x\in M_{\Omega}(n)$ then $\pi_n(z) = x$ if and only if $$(\rho_{1\,m}(z),\ldots,
\rho_{m\,m}(z)) = (\mu_{1\,m}(x),\ldots, \mu_{m\,m}(x))\eqno (4.10)$$ up to a reordering, for
all $m\in I_n$, using the notation of \S 2.2. \vs One defines a free action of $\Sigma_n$ on
$M_{\Omega}(n,\ee)$, operating as a group of algebraic isomorphisms, by defining $$\sigma\cdot
z = (\sigma\cdot\nu, x)\eqno (4.11)$$ where $\sigma\in \Sigma_n $ and $z=(\nu,x)\in
M_{\Omega}(n,\ee)$. \vs The following well--known proposition is classical.\vs {\bf Proposition
4.2.} {\it Let $m\in I_n$ and let $\gamma(m)\in \ee(m)$. Let $g(m)\in Gl(m)$ and put
$x(m)=Ad\,g(m)(\gamma(m))$. Then there exists an open neighborhood $N$ of $\gamma(m)$ in
$\ee(m)$ and a section $S\s Gl(m)$ of the quotient map $Gl(m)\to Gl(m)/Diag(m)$ (using
notation in
\S3.4) defined on a neighborhood of $g(m)\,Diag(m)$ such that the map $$S\times N\to
M(m)\,\,\hbox{where}\,\,(g,\gamma)\mapsto Ad\,g(\gamma)\eqno (4.12)$$ is an analytic isomorphism onto
an open set of $M(m)$. The elements of the image are necessarily regular semisimple elements of
$M(m)$.}\vs If $U$ is an open subset of $M_n(\Omega)$ and $\phi:U\to \ee$ is an analytic map
let $graph\,\phi: U\to \ee\times M_n(\Omega)$ be the analytic map defined by putting
$graph\,\phi(x) = (\phi(x),x)$. \vs
{\bf Proposition 4.3.} {\it Let $z = (\nu,x)\in M_{\Omega}(n,\ee)$. Then there exists a
(sufficiently small) connected open neighborhood $U$ of $x$ in $M_{\Omega}(n)$ and an analytic
map
$\phi:U\to
\ee_{\Omega(n)}$ with the following properties:

(1)$\,\,graph\,\phi(x) = z,\,\,graph\,\phi:U\to U_z$ is a homeomorphism, where $U_z$
is the image of $graph\,\phi$, and $U_z\s M_{\Omega}(n,\ee)$.

(2)$\,\,\pi_n^{-1}(U) =\sqcup_{\sigma\in \Sigma_n}\sigma\cdot U_z$. 

(3)$\,\,\{\sigma\cdot
U_z\mid \sigma\in \Sigma_n\}$ are the connected components of $\pi_n^{-1}(U)$ and each component is
open in
$M_{\Omega}(n,\ee)$.

In particular $U$ is evenly covered by $\pi_n$ and $\pi_n$ is a covering projection (see
(4.8)).}\vs {\bf Proof.} Statements (1) and (2) are immediate consequences of Proposition 4.2. But it
is immediate from (1) and (2) that $\sigma\cdot U_z$ is connected and closed in $\pi_n^{-1}(U)$ for
any $\sigma\in \Sigma_n$. Since the partition in (2) is finite it follows that the parts are open in 
$\pi_n^{-1}(U)$ and hence are open in $M_{\Omega}(n,\ee)$. The remaining statements are obvious. QED 
\vs 4.3. As in the introduction, \S 0, let $\b_e = -e +\b$ using the notation of \S2.2. By
Remark 2.4, Theorems 2.3 and 2.5 hold if $\b_e$ replaces $e + \b$. Let $\b_{e,\Omega(n)} =
M_{\Omega}(n)\cap \b_e$. Then $\b_{e,\Omega(n)}$ is a Zariski open subset of $\b_e$ and
$$\Phi_n:\b_{e,\Omega(n)}\to \Omega(n)\eqno (4.13)$$ is an algebraic isomorphism by Theorem
2.5 (see Remark 2.16). In particular $\b_{e,\Omega(n)}$ is dense in $\b_e$. Now let
$M_{\Omega}(n,\ee,\b) =
\pi_n^{-1}(\b_{e,\Omega(n)})$ so that, by Proposition 4.3, $M_{\Omega}(n,\ee,\b)$ is a covering
of $\b_{e,\Omega(n)}$. Now consider the restriction $$\kappa_n:M_{\Omega}(n,\ee,\b)\to 
\ee_{\Omega(n)}\eqno (4.14)$$ of (4.9) to $M_{\Omega}(n,\ee,\b)$. Note that $\ee_{\Omega(n)}$
is connected since it is clearly Zariski open in $\ee$. One has
\vs {\bf Theorem 4.4.} {\it The map (4.14) is a homeomorphism. In particular
$M_{\Omega}(n,\ee,\b)$ is connected.}\vs {\bf Proof.} Let $\nu\in \ee_{\Omega(n)}$. Then by
Theorem 2.5 there exists (uniquely) $x\in \b_e$ such that, for any $m\in I_n,\,\,
(\mu_{1\,m}(x),\ldots,\mu_{m,m}(x)) = (\nu_{1,\,m},\ldots,\nu_{m\,m})$, up to a reordering.
But then $x\in \b_{e,\Omega(n)}$ and, by Proposition 4.3, $(\nu,x)\in M_{\Omega}(n,\ee,\b)$.
But then $\nu$ is in the image (4.14). That is, (4.14) is surjective. But assume $z,z'\in
M_{\Omega}(n,\ee,\b)$ and $\kappa_n(z) = \kappa_n(z')$. Then if $x = \pi_n(z)$ and 
$x'=\pi_n(z')$ one has $x,x'\in \b_e$ and hence $x=x'$ by Theorem 2.5. Thus $z=z'$ so that
(4.14) is injective. Hence (4.14) is bijective. But of course as a restriction map, (4.14) is
continuous. We have only to show that its inverse is continuous. 

Let $\beta_i,\,i\in I_n$, be the regular function on $\ee_{\Omega(n)}$ defined so that if
$\nu\in \ee_{\Omega(n)}$ then $\beta_{d(m-1) + k}(\nu),\,\,k\in I_m,\,m\in I_n$, is the
elementary symmetric function of degree $m-k +1$ in $\{\nu_{1\,m},\ldots,\nu_{m\,m}\}$. Now
let 
$\beta:\ee_{\Omega(n)}\to \Bbb C^{d(n)}$ be the regular algebraic map defined so that
$$\beta(\nu) = (\beta_1(\nu),\ldots, \beta_{d(n)}(\nu))\eqno (4.15)$$ One notes that if
$c=\beta(\nu)$ then by (2.3),(2,4),(2.10) and (2.11) $$(\nu_{1\,m},\ldots,\nu_{m\,m}) =
(\mu_{1\,m}(c),\,\ldots,\mu_{m\,m}(c))\eqno (4.16)$$ up to a reordering, for all $m\in I_n$.
It follows then that $$\beta:\ee_{\Omega(n)}\to \Omega(n)\eqno (4.17)$$ is a surjective
morphism (see (2.53)). Recalling Theorem 2.3 (where $-e$ replaces $e$), one has, inverting
(4.13), a surjective morphism $$\widetilde {\beta}:\ee_{\Omega(n)}\to
\b_{e,\Omega(n)}\eqno (4.18)$$ where $\beta = \Phi_n\circ \widetilde {\beta}$ noting that,
for $i\in I_n$,
$$p_i(\widetilde {\beta}(\nu)) = \beta_i(\nu)\eqno (4.19)$$ by (2.3),(2.4) and (2.5). But
clearly $(\nu,\widetilde{\beta}(\nu))\in M_{\Omega}(n,\ee,\b)$ for any $\nu\in
\ee_{\Omega(n)}$. Hence $$\ee_{\Omega(n)}\to M_{\Omega}(n,\ee,\b),\,\,\,\nu\mapsto
(\nu,\widetilde{\beta}(\nu))\eqno (4.20)$$ is an algebraic morphism. But (4.20) must be the inverse
to (4.14), by the bijectivity of (4.14), since
$\kappa_n((\nu,\widetilde{\beta}(\nu))) = \nu$. Hence (4.14) is a homeomorphism. QED\vs Recalling
(4.9) let $M_{\nu}(n,\ee) = \kappa_{n}^{-1}(\nu)$ so that one has a ``fibration" 
$$M_{\Omega}(n,\ee) =\sqcup_{\nu\in \ee_{\Omega (n)}} M_{\nu}(n,\ee)\eqno (4.21)$$ of
$M_{\Omega}(n,\ee)$ over $\ee_{\Omega(n)}$ with fiber projection $\kappa_n$ (see (4.9) and
Proposition 4.5 below). Let
$\nu\in
\ee_{\Omega(n)}$. If $\nu \in \ee_{\Omega(n)}$ and $c\in \Omega(n)\s \Bbb C^{d(n)}$ is defined
by
$$c=
\beta(\nu)\eqno (4.22)$$ (see (4.17)) note that $$M_c(n) \to M_{\nu}(n,\ee),\qquad x\mapsto
 (\nu,x)\eqno (4.23)$$
is a homeomorphism, by the definition of
$M_{\Omega}(n,\ee)$. The following asserts, in particular, that the ``fibers" of (4.21) are all
homeomorphic. \vs {\bf Proposition 4.5.} {\it  One has a
homeomorphism $$ M_{\nu}(n,\ee)\cong (\Bbb C^{\times})^{d(n-1)}\eqno (4.24)$$ for any
$\nu\in\ee_{\Omega(n)}$.}\vs {\bf Proof.} This is immediate from (4.23) and Theorem 3.23. QED
\vs {\bf Remark 4.6.} Note that $M_{\Omega}(n,\ee,\b)$ defines a
cross--section of $\kappa_n$ by Theorem 4.14. That is, for any $\nu\in
\ee_{\Omega(n)}$ the intersection
$$M_{\Omega}(n,\ee,\b)\cap M_{\nu}(n,\ee)\,\,\hbox{has only one point}\eqno (4.25)$$\vskip
.5pc
{\bf Proposition 4.7.}
{\it $M_{\Omega}(n,\ee)$, as a topological space (Euclidean topology), is connected. }\vs 
{\bf Proof.} As a covering of the manifold of $M_{\Omega}(n)$ obviously $M_{\Omega}(n,\ee)$ is
locally connected (see Proposition 4.3) so that any connected component of
$M_{\Omega}(n,\ee)$ is open in $M_{\Omega}(n,\ee)$. But then there exists a connected
component $C$ such that
$M_{\Omega}(n,\ee,\b)\s C$, by Theorem 4.4. But if $\nu \in \ee_{\Omega(n)}$
then $M_{\nu}(n,\ee)$ is connected by (4.24). But then $M_{\nu}(n,\ee)\s C$ by (4.25). Hence
$C= M_{\Omega}(n,\ee)$ by (4.21). Thus $M_{\Omega}(n,\ee)$ is connected. QED\vs 4.4. The definition
 of variety here
and throughout implies that it is Zariski irreducible. We will
prove in this section that $M_{\Omega}(n,\ee)$ is a nonsingular affine variety of dimension
$n^2$. We first observe \vs {\bf Proposition 4.8.} {\it The nonempty Zariski open subset
$\ee_{\Omega(n)}$ of $\ee$ (see
\S 4.1) is a nonsingular affine variety and the nonempty Zariski open subset $M_{\Omega}(n)$ of
$M(n)$ is a nonsingular affine variety (see (2.53)). In particular $\ee_{\Omega(n)}\times
M_{\Omega}(n)$ is a nonsingular affine variety of dimension $n^2 + d(n)$.}\vs {\bf Proof.} Recalling
\S 4.1, clearly, for any $m\in I_n$, $\ee(m)$ is a Zariski open, nonempty subvariety of $\d(m)$. It
is affine since it is the complement of the zero set of the discriminant function on $\d(m)$. But then
$\ee$ is a nonsingular affine variety of dimension $d(n)$. But then $\ee_{\Omega(n)}$ is a nonsingular
affine variety of dimension $d(n)$ since the condition (4.4) clearly defines $\ee_{\Omega(n)}$
as the complement of the zero set of a single regular function on $\ee$. But now the
argument in Remark 2.16 readily characterizes $M_{\Omega}(n)$ as the complement in $M(n)$ of the zero
set in 
$M(n)$ of a polynomial in $J(n)$ (see (2.30)). Thus $M_{\Omega}(n)$ is a nonsingular affine variety of
dimension $n^2$. This of course proves the proposition. QED\vs  If $X$ is an affine variety we will denote
the affine ring of $X$ by ${\cal O}(X)$. \vs {\bf Theorem 4.9.} {\it
$M_{\Omega}(n,\ee)$ is a $n^2$-dimensional Zariski closed affine nonsingular subvariety of the
nonsingular affine variety
$\ee_{\Omega(n)}\times M_{\Omega}(n)$ (see Proposition 4.8).}\vs {\bf Proof.} One has 
$${\cal O}(\ee_{\Omega(n)}\times M_{\Omega}(n)) = {\cal O}(\ee_{\Omega(n)}) \otimes
{\cal O}(M_{\Omega}(n))\eqno (4.26)$$ For $i\in I_n$ let $\beta_i\in {\cal O}(\ee_{\Omega(n)})$  be
defined as in  (4.15). Also let $p_i'
\in {\cal O}(M_{\Omega}(n))$ be defined by putting $p_i' = p_i|M_{\Omega}(n)$ where $p_i\in J(n)$ is
given by (2.5) (see (2.30)). For notational convenience put $W = \ee_{\Omega(n)}\times M_{\Omega}(n)$.
Now let
$f_i\in {\cal O}(W)$ be defined by putting $f_i = \beta_i\otimes 1 - 1\otimes p_i'$. Clearly
$$M_{\Omega}(n,\ee) = Spec[{\cal O}(W)/(f_1,\ldots,f_{d(n)}]\eqno (4.27)$$ But for any $z=(\nu,x)\in
M_{\Omega}(n,\ee)$ one has $$(df_i)_z,\,i\in I_n,\,\,\,\hbox{are linearly independent}\eqno (4.28)$$
since $(dp_i)_x,\,i\in I_n$, are linearly independent in $T^*_x(M(n))$ by (2.55) and the definition of
$M^{sreg}(n)$ in
\S 2.3. One also notes that $(d\beta_i)_{\nu},\,i\in I_n$, are linearly independent in $T^*_{\nu}(\ee)$
since $\nu(m)$ is regular in $\d(m)$ for any $m\in I_n$. But then $M_{\Omega}(n,\ee)$ is nonsingular
at $z$ by Theorem 4 in
\S 4 of Chapter 3 in [M1], p. 172, where $X$ and $U$ in the notation of that reference are equal to $W$
here (see Proposition 4.8) and $Y = M_{\Omega}(n,\ee)$. Thus $M_{\Omega}(n,\ee)$ is nonsingular and
has dimension $n^2$. But now since $M_{\Omega}(n,\ee)$ is connected in the Euclidean topology, by
Proposition 4.7, it is obviously connected in the Zariski topology (i.e. it is not the disjoint union
of two nonempty Zariski open sets).  But then, since it is nonsingular, it is irreducible as an
algebraic set by Corollary 17.2, p. 72, in [B]. It is then also a closed affine subvariety of $W$
by (4.27). QED\vs We may regard ${\cal O}(M_{\Omega}(n))$ as a module for $J(n)$ (see \S 2.4)) where
$p_i$ operates as multiplication by $p_i'$ (using the notation in the proof of Theorem 4.9).
Similarly regard ${\cal O}(\ee_{\Omega(n)})$ as a module for $J(n)$ where $p_i$ operates as
multiplication by $\beta_i$ (see (4.19)). Then Theorem 4.9 and the equality (4.27) clearly implies \vs
{\bf Theorem 4.10.} {\it As an affine ring one has the tensor product $${\cal O}(M_{\Omega}(n,\ee)) = 
{\cal O}(\ee_{\Omega(n)}) \otimes_{J(n)}{\cal O}(M_{\Omega}(n))\eqno (4.29)$$}\vs 4.5. 
Assume $Y$ is an affine variety and that $\Sigma$ is a finite group
operating as a group of algebraic isomorphisms of $Y$. Then $\Sigma$ operates as a group of algebraic
automorphisms of ${\cal O}(Y)$ where for any $f\in {\cal O}(Y),\,\sigma\in \Sigma $ and $z\in Y$,
one has $$(\sigma\cdot f)(z) = f(\sigma^{-1}\cdot z)\eqno (4.30)$$ Let $Y/\Sigma$ be the set of
orbits of $\Sigma$ on $Y$. Obviously any $f\in {\cal O}(Y)^{\Sigma}$ defines a function of $Y/\Sigma$.
It is then a classical theorem that $Y/\Sigma$ has the structure of an affine variety where $${\cal
O}(Y/\Sigma) = {\cal O}(Y)^{\Sigma}\eqno (4.31)$$ See e.g. Theorem 1.1, p. 27 and Amplification 1.3,
p. 30 in  Chapter 1,
\S 2 of [M2]. In
addition one notes that $$Y\to Y/\Sigma,\,\,z \mapsto \Sigma\cdot z \eqno (4.32)$$ is a morphism where
the corresponding cohomorphism is the embedding $${\cal O}(Y)^{\Sigma}\to {\cal O}(Y)\eqno (4.33)$$
Let $F(Y)$ be the quotient field of ${\cal O}(Y)$ and let $F(Y/\Sigma)$ be the quotient field of
${\cal O}(Y/\Sigma)$. Since $\Sigma$ is finite it is a simple and well-known fact that, as a consequence
of (4.31), $$F(Y/\Sigma) = F(Y)^{\Sigma}\,\,\,\hbox{so that $F(Y)$ is a Galois extension of
$F(Y/\Sigma)$}\eqno (4.34)$$ Let $Clos_{F(Y)}({\cal O}(Y/\Sigma)$ be the integral closure of ${\cal
O}(Y/\Sigma)$ in $F(Y)$. \vs {\bf Proposition 4.11.} {\it $Clos_{F(Y)}({\cal O}(Y/\Sigma)$ is a finite
module over ${\cal O}(Y/\Sigma)$. In particular $Clos_{F(Y)}({\cal O}(Y/\Sigma)$ is Noetherian.
Furthermore ${\cal O}(Y)$ is also a finite module over ${\cal O}(Y/\Sigma)$ so that ${\cal O}(Y)$ is
integral over ${\cal O}(Y/\Sigma)$ and hence $${\cal O}(Y)\s Clos_{F(Y)}({\cal O}(Y/\Sigma)\eqno
(4.35)$$ and one has equality in (4.35) in case $Y$ is nonsingular. Finally, (in the sense of Chapter 2,
\S 7, Definition 3, p. 124 in [M1]) the
morphism (4.32) is finite and the morphism
$$Z\mapsto Y/\Sigma\eqno (4.36)$$ is finite where $Z= Spec(Clos_{F(Y)}({\cal O}(Y/\Sigma)))$ and
(4.36) is defined so that the injection ${\cal O}(Y/\Sigma)\to Clos_{F(Y)}({\cal O}(Y/\Sigma)$ is the
corresponding cohomomorphism.}\vs {\bf Proof. } The first statement is given by Theorem 9, in Chapter
5, \S 4, p.267 in [ZS]. The
statement that ${\cal O}(Y)$ is also a finite module over ${\cal O}(Y/\Sigma)$ is stated as
Noether's Theorem and proved as Theorem 2.3.1, p. 26 in [Sm]. But now if $Y$ is nonsingular then
$O(Y)$ is integrally closed in $F(Y)$. See e.g. the first paragraph, p. 197, in [M1]. But of course 
$Clos_{F(Y)}({\cal O}(Y/\Sigma)$ is integral over $O(Y)$. Hence one has equality in (4.35). But now
the finiteness of (4.32) and (4.36) follows from Proposition 5, p. 124 in \S 7 of Chapter 2 in [M1]
since $Y, Y/\Sigma$ and $Z$ are affine varieties. QED\vs Making use of Theorem 4.9 we apply Proposition
4.11 in the case where
$Y = M_{\Omega}(n,\ee)$ and $\Sigma= \Sigma_n$ with, of course, the action given by (4.12). The
quotient field of ${\cal O}(M_{\Omega}(n,\ee))$ will be denoted by $F(n,\ee)$. The quotient field 
of ${\cal O}(M(n))$ will be denoted by $F(n)$. In the
notation of \S 1.1 note that
${\cal O}(M(n))$ is just $P(n)$. Since $M_{\Omega}(n)$ is
Zariski dense in $M(n)$ note that $F(n)$ is also the quotient field of ${\cal O}(M_{\Omega}(n))$.

It is clear, from the definition of
$\pi_n$ (see (4.8)), that 
$\pi_n$ is
a morphism whose corresponding cohomomorphism maps ${\cal O}(M_{\Omega}(n))$ into ${\cal
O}(M_{\Omega}(n,\ee)^{\Sigma_n})$. Hence $\pi_n$ descends to a morphism
$$M_{\Omega}(n,\ee)^{\Sigma_n} \to M_{\Omega}(n)\eqno (4.37)$$ \vs {\bf Proposition 4.12.} {\it The
morphism (4.37) is an isomorphism of algebraic varieties. That is, $$M_{\Omega}(n,\ee)^{\Sigma_n}
\cong  M_{\Omega}(n)\eqno (4.38)$$ In particular $\pi_n$ is a finite morphism. Furthermore
$$F(n,\ee)^{\Sigma_n}\cong F(n) \eqno (4.39)$$ so that, using (4.39) to define an identication,
$F(n,\ee)$ is a Galois extension of $F(n)$ with Galois group $\Sigma_n$. In addition, using (4.38) to
define an identifation one has $$Clos_{F(n,\ee)}{\cal O}(M_{\Omega}(n)) = {\cal
O}(M_{\Omega}(n,\ee))\eqno (4.40)$$}\vs {\bf Proof.} It is immediate from Proposition 4.3 that (4.37)
is bijective. But then it is birational by Theorem 5.1.6, p. 81 in [Sp] (since we are in a
characteristic zero case). But
$M_{\Omega}(n)$ is nonsingular. Thus (4.37) is an isomorphism (see e.g. Theorem 5.2.8, p. 85 in [Sp]).
The rest of the statements follow from Proposition 4.11 since $M_{\Omega}(n,\ee)$ is nonsingular by
Theorem 4.9. QED\vs Recall the notation of Proposition 4.3 so that $z= (\gamma,x)\in
M_{\Omega}(n,\ee)$. Also $U_z$ is a (Euclidean) open neighborhood of $z$ in $M_{\Omega}(n,\ee)$, $U$
is a (Euclidean) open neighborhood of $x$ in $M_{\Omega}(n)$ and the statements of Proposition 4.3
hold. The inverse of the homeomorphism $graph\,\phi:U\to U_z$ is clearly $$\pi_n|U_z:U_z\to U\eqno
(4.41)
$$
\vskip.5pc {\bf Proposition 4.13.} {\it Recalling the notation of Proposition 4.3 the map (4.41) is an
analytic isomorphism. In particular (see Corollary 2, p. 182, to Theorem 3 in \S 5 of Chapter 3 in 
[M1]) $\pi_n: M_{\Omega}(n,\ee) \to M_{\Omega}(n)$ (see (4.8)) is an \'etale morphism.}\vs 
{\bf Proof.}
Since $\pi_n$ is a morphism it is a holomorphic map of nonsingular analytic manifolds (see ii, p. 58
of \S 10, Chapter 1 in [M1]) Thus the homeomorphism (4.41) is an analytic map. It suffices to prove that
$$graph\,\phi:U\to U_z\eqno (4.42)$$ is analytic. To do this first regard (4.42) as a map $$ U\to
\ee_{\Omega(n)}\times M_{\Omega}(n)\eqno (4.43)$$ (see (4.7)). By the definition of $graph\,\phi$ in
Proposition 4.3 it is obvious that (4.43) is analytic. Hence if $g\in {\cal O}(\ee_{\Omega(n)}\times
M_{\Omega}(n))$ then $g\circ graph\,\phi$ is analytic function on $U$. But then $f\circ 
graph\,\phi$ is an analytic function on $U$ for any $f\in {\cal O}(M_{\Omega}(n,\ee))$ since ${\cal
O}(M_{\Omega}(n,\ee))$ is just the restriction of ${\cal O}(\ee_{\Omega(n)}\times M_{\Omega}(n))$ to
$M_{\Omega}(n,\ee)$. But then (4.42) is analytic since an analytic coordinate system in a
Euclidean neighborhood of
$z$ in
$M_{\Omega}(n,\ee)$ is given by elements in ${\cal O}(M_{\Omega}(n,\ee))$ which are uniformizing
parameters in a Zariski neighborhood of $z$ (see p. 183, \S 6 of Chapter 3 in [M1]). QED\vs
Combining Propositions 4.3 and 4.13 one has \vs {\bf Proposition 4.14.} {\it The map
$\pi_n$ defines $M_{\Omega}(n,\ee)$ an analytic covering of $M_{\Omega}(n)$ with $\Sigma_n$ as
the group of deck transformations.}\vs 4.6. Since
$M_{\Omega}(n,\ee)$ is locally and analytically isomorphic to $M_{\Omega}(n)$ (via $\pi_{\n}$) the
tensor which defines Poisson bracket of functions on $M_{\Omega}(n)$ lifts and defines Poisson bracket
of analytic functions on $M_{\Omega}(n,\ee)$. In particular ${\cal O}(M_{\Omega}(n,\ee))$ has the
structure of a Poisson algebra. For any $f\in {\cal O}(M_{\Omega}(n))$ (noting that we regard ${\cal
O}(M(n))\s {\cal O}(M_{\Omega}(n))$) let
$\widehat{f}\in {\cal O}(M_{\Omega}(n,\ee))$ be defined by putting $\widehat{f} = f\circ \pi_n$. For
$f_1,f_2 \in {\cal O}(M_{\Omega}(n))$ one then has $$\widehat {[f_1,f_2]} = [\widehat {f_1},\widehat
{f_2}]\eqno (4.44)$$ Using the notation of \S 1.2 but now, in addition, applied to $M_{\Omega}(n,\ee)$,
for any
$\varphi \in {\cal O}(M_{\Omega}(n,\ee))$ let $\xi_{\varphi}$ be the (complex) vector field on
$M_{\Omega}(n,\ee)$ defined so that $\xi_{\varphi}\psi = [\varphi,\psi]$ for any $\psi\in {\cal
O}(M_{\Omega}(n))$. It is immediate that if $f\in {\cal O}(M_{\Omega}(n)$ then $\xi_{\widehat f}$ is
$\pi_n$-related to $\xi_f$ so that unambiguously $$(\pi_n)_*(\xi_{\widehat f}) = \xi_f\eqno (4.45)$$
Besides the (just considered) subring of ${\cal O}(M_{\Omega}(n,\ee))$, defined by the pull-back of (the
surjection) $\pi_n$, there is the subring of ${\cal O}(M_{\Omega}(n,\ee))$ defined by the pull-back of
(the surjection)
$\kappa_n$ (see (4.9)). Indeed let $$J(n,\ee) = \{q\circ \kappa_n\mid q\in {\cal
O}(\ee_{\Omega(n)})\}\eqno (4.46)$$ From the definition of $\rho_{k\,m},\,k\in I_m,\,\,m\in I_m$ in \S
4.2 note that $\rho_{k\,m} \in J(n,\ee)$. For notational convenience let $r_i\in J(n,\ee),\,\,i\in
I_{d(n)}$ be defined so that $$r_i = \rho_{k\,m}\eqno (4.47)$$ where $$i = d(m-1) + k\eqno (4.48)$$ so
that $$r_i\in J(n,\ee),\,\,i\in d(n)\eqno (4.49)$$\vskip .5pc {\bf Remark 4.15.} Recalling \S 4.1
note that conversely $J(n,\ee)$ is a localization of the polynomial ring generated by the 
$r_i,\,i\in I_{d(n)}$. \vs Now put $\widehat {J(n)} = \{\widehat {p}\mid p\in J(n)\}$ (see \S  2.4) so
that $\widehat {J(n)}$ is the polynomial ring $$\widehat {J(n)} = \Bbb C[\widehat {p_1},\ldots,
\widehat {p_{d(n)}}]\eqno (4.50)$$ Let $I_{[m]} =
I_{d(m)}-I_{d(m-1)}$ so that $card\,\,I_{[m]} = m.$ \vs {\bf Proposition 4.16.} {\it One has 
$$\widehat
{J(n)}\s J(n,\ee)\eqno (4.51)$$ In fact if $i\in I_{[m]}$, where $m\in I_n$, and $i$ is written as in
(4.48) then $\widehat {p_i}$ is the elementary symmetric polynomial of degree $m-k+1$ in
the functions
$r_j,\,j\in I_{[m]}$. Indeed if $$P_{m}(\lambda) = \lambda^m + \sum_{k\in I_{m}}
(-1)^{m-k+1}\widehat {p_{d(m-1) +k}}\,\,\lambda^{k-1}$$ then $$ P_{m}(\lambda) = \prod_{j\in
I_{[m]}}(\lambda- r_j)\eqno (4.52)$$ so that, in addition, $r_j$, for
$j\in I_{[m]}$, satisfies the polynomial equation
$$P_m(r_j)=0 $$}\vs {\bf Proof.} The inclusion (4.51) follows from (4.50) and (4.52). On the
other hand (4.52) follows from (2.3), (2.4) and (2.5) together with (4.10), (4.29) and (4.47). QED\vs 
Since $\pi_n$ is an analytic covering map (see Proposition 4.14) it follows from (2.55) and the
definition of strongly regular (see \S 2.3) that the differentials $(d\widehat {p_i})_z,\,i\in
I_{d(n)}$, are linearly independent at any $z\in M_{\Omega}(n,\ee)$. For any $m\in I_n$ let
$T_z^*(M_{\Omega}(n,\ee))^{(m)}$ be the $m$-dimensional subspace of the cotangent space
$T_z^*(M_{\Omega}(n,\ee))$ spanned by the differentials $(d\widehat {p_i})_z,\,i\in
I_{[m]}$. 
\vs {\bf Proposition 4.17.} {\it Let $z\in M_{\Omega}(n,\ee)$ and let $m\in I_n$. 
Then $(dr_i)_z,\,i\in
I_{[m]}$ is a basis of $T_z^*(M_{\Omega}(n,\ee))^{(m)}$. }\vs {\bf Proof.} Let $V$ be the space of 
$T_z^*(M_{\Omega}(n,\ee))$ spanned by $(dr_i)_z,\,i\in I_{[m]}$. But then $dim\,V \leq m$. But 
$T_z^*(M_{\Omega}(n,\ee))^{(m)}\s V$ by (4.51). Thus $V= T_z^*(M_{\Omega}(n,\ee))^{(m)}$ by dimension.
QED\vs In the notation above let $T_z^*(M_{\Omega}(n,\ee))'$ (resp. $T_z^*(M_{\Omega}(n,\ee))''$) be 
the $d(n-1)$-(resp $d(n)$-) dimensional sum of the subspaces $T_z^*(M_{\Omega}(n,\ee))^{(m)}$ over all
$m\in I_{n-1}$ (resp. $m\in I_n$). Then as an immediate consequence of Proposition 4.17 one has \vs
{\bf 
Proposition 4.18.} {\it Let $z\in M_{\Omega}(n,\ee)$. Then $(dr_i)_z,\,i\in
I_{d(n-1)}$ (resp. $I_{d(n)}$) is a basis of $T_z^*(M_{\Omega}(n,\ee))'$ (resp.
$T_z^*(M_{\Omega}(n,\ee))'')$.} \vs 4.7.  Let
$\nu \in \ee_{\Omega(n)}$. By definition (see (4.21)) $M_{\nu}(n,\ee) = \kappa_n^{-1}(\nu)$ so that 
$M_{\nu}(n,\ee)$ is a Zariski closed subset of $M_{\Omega}(n,\ee)$. \vs {\bf Remark 4.19.} Note that
(see (4.2), (4.3), \S 4.2 and (4.47))
$M_{\nu}(n,\ee)$ may be given by the equations $$M_{\nu}(n,\ee) = \{z\in M_{\Omega}(n,\ee) \mid r_i(z)
= \nu_{k\,m}\,\hbox{when $i\in I_{d(n)}$ is put in the form (4.48)}\}\eqno (4.53)$$ Note also
that $M_{\nu}(n,\ee)$ is nonsingular by Proposition 4.18 (linear independence of differentials). \vskip
1pc {\bf Proposition 4.20.} {\it Let
$\nu
\in
\ee_{\Omega(n)}$ and let
$c=
\beta(\nu)$ so that $c\in
\Omega(n)$ (see (4.22)). Then the covering map $\pi_n$ (see (4.8)) restricts to an algebraic
isomorphism $$\pi_n:  M_{\nu}(n,\ee) \to M_c(n)\eqno (4.54)$$ of nonsingular affine varieties. }\vs
{\bf Proof.} Recall that $M_c(n)$ is an irreducible nonsingular Zariski closed subvariety of $M(n)$ (see
Theorem 3.23). The homeomorphism (4.23) can obviously be regarded as a morphism mapping $M_c(n)$ to
$\ee_{\Omega(n)}\times M_{\Omega}(n)$. However the image of (4.23) is the Zariski closed subset
$M_{\nu}(n,\ee)$ of 
$\ee_{\Omega(n)}\times M_{\Omega}(n)$. Thus (4.23), as it stands, is a bijective morphism. But
then 
$M_{\nu}(n,\ee)$ is a variety (i.e. it is irreducible). Hence it is a nonsingular affine variety by
Remark 4.19. Thus (4.23), as it stands, is an
algebraic isomorphism. But (4.54) is just the inverse of (4.23). QED\vs Note that (2.31) and 
(4.44) imply $$[\widehat {p},\widehat {q}] = 0\eqno (4.55)$$ for any $p,q\in J(n)$. In particular 
$$[\widehat {p_i},\widehat {p_j}] = 0\eqno (4.56)$$ for any $i,j\in I_{d(n)}$. One consequence of
Proposition 4.19 is \vs {\bf Proposition 4.21.} {\it  Let $\nu \in \ee_{\Omega(n)}$ and let 
$z\in M_{\nu}(n,\ee)$, so that $\nu =\kappa_n(z)$ (see (4.9)). Then $(\xi_{\widehat p_i})_z,\,\,i\in
I_{d(n-1)}$, is a basis of the tangent space $T_z(M_{\nu}(n,\ee))$.} \vs {\bf Proof.} Let $x =
\pi_n(z)$ (see (4.8)) and $c = \beta(\nu)$ (see (4.22)) so that $x\in M_c(n)$. Since $\pi_n$ is a
local analytic isomorphism it suffices by Proposition 4.20 and (4.45) to see that
$(\xi_{p_i})_x,\,i\in d(n-1)$, is a basis of $T_x(M_c(n))$. But $x$ is strongly regular (see \S 2.3)
by (2.55) since $c\in \Omega(n)$. The result then follows from Remark 2.8 and Theorems 3.4, 3.23.
QED\vs The argument which established Theorem 3.25 may now be used to establish \vs {\bf Theorem
4.22.} {\it $J(n,\ee)$ (see (4.46) is a maximal Poisson commutative subalgebra of ${\cal
O}(M_{\Omega}(n,\ee))$. In particular (see (4.48)) $$[r_i,r_j] = 0\eqno (4.57)$$ for any $i,j\in
I_{d(n)}$. Furthermore if $f\in J(n,\ee)$ and $\nu\in \ee_{\Omega(n)}$ then $f|M_{\nu}(n,\ee)$ is a
constant function and $\xi_f|M_{\nu}(n,\ee)$ is tangent to $M_{\nu}(n,\ee)$. Moreover if we write
(using Proposition 4.21)
$$\xi_f  = \sum_{i\in I_{d(n-1)}}\, f_i\,\,\xi_{\widehat {p_i}}\eqno (4.58)$$ on $M_{\nu}(n,\ee)$ where
$f_i\in {\cal O}(M_{\nu}(n,\ee))$ then all the $f_i$ are constant on $M_{\nu}(n,\ee)$. Finally
$$\eqalign{(\xi_{r_i})_z&,\,i\in I_{d(n-1)}, \,\,\hbox{is a basis of $T_z(M_{\nu}(n,\ee))$ for any
$z\in M_{\nu}(n,\ee)$ and}\cr
(\xi_{r_i})_z&=0,\,i\in I_{[n]}= I_{d(n)}-I_{d(n-1)},\,\,\hbox{for any
$z\in M_{\nu}(n,\ee)$}\cr}\eqno (4.59)$$}\vs {\bf Proof.} Let $\nu \in \ee_{\Omega(n)}$. The function
$r_i$ is constant on
$M_{\nu}(n,\ee)$ for all $i\in I_{d(n)}$ by Remark 4.19. Let $f\in J(n,\ee)$. But then
$f|M_{\nu}(n,\ee)$ is a constant function by Remark 4.15. On the other hand if $z\in  M_{\nu}(n,\ee)$
and $W_z$ is the orthocomplement of $T_z(M_{\nu}(n,\ee))$ in $T^*(M_{\Omega}(n,\ee))$ then
$(dr_i)_z,\,i\in I_{d(n)}$, is a basis of $W_z$ by Proposition 4.18 and Remark 4.19. But this implies
that $(df)_z\in W_z$ by Remark 4.15. In particular $(d\widehat {p_i})_z\in W_z$ for any $i\in I_{d(n)}$
by (4.50). In fact
$(d\widehat {p_i})_z,\,i\in I_{d(n)}$, is a basis of $W_z$ by Proposition 4.18 and the definition of
$T_z^*(M_{\Omega}(n,\ee))''$ in \S 4.6. 

Now for any $g\in
{\cal O}(M_{\Omega}(n,\ee)$ the tangent vector $(\xi_g)_z$ depends only on $(dg)_z$ (see \S 1.2). But
 $$ (\xi_{\widehat {p_i}})_z = 0\,\,\hbox{for any $i\in I_{[n]}= I_{d(n)}-I_{d(n-1)}$}\eqno (4.60)$$
by (2.7). But then $(\xi_g)_z\in T_z(M_{\nu}(n,\ee))$ if $(dg)_z\in W_z$, by Proposition 4.21. Hence
$\xi_f|M_{\nu}(n,\ee)$ is tangent to $M_{\nu}(n,\ee)$. Furthermore Propositions 4.17,4.18 and 4.21 
imply (4.59). Now if $g\in J(n,\ee)$
then $g|M_{\nu}(n,\ee)$ is constant. But $\xi_f|M_{\nu}(n,\ee)$ is tangent to $M_{\nu}(n,\ee)$. Thus 
$J(n,\ee)$ is Poisson commutative. In particular $[\xi_{\widehat {p_j}},\xi_{f}] = 0$ for any $j\in
I_{d(n-1)}$. But, on $M_{\nu}(n,\ee)$, one has $$[\xi_{\widehat {p_j}},\xi_{f}] = \sum_{i\in
I_{d(n-1)}}\, \,(\xi_{\widehat {p_j}}f_i)\,\,\xi_{\widehat {p_i}} $$ by
(4.58). Thus $\xi_{\widehat {p_j}}f_i = 0$, by Proposition 4.21, for all $i,j\in I_{d(n-1)}$. Hence the
$f_i$ are constants. 

Now recall the definition of $M_{\Omega}(n,\ee,\b)$ in \S 4.3. Then $M_{\Omega}(n,\ee,\b)$ is a Zariski
closed subset of $M_{\Omega}(n,\ee)$ since, clearly, $\b_{\ee,\Omega(n)}$ is obviously closed in
$M_{\Omega}(n)$. But $M_{\Omega}(n,\ee,\b)$ is irreducible since it is the image of the  bijective
morphism (4.20). But then (4.14) is a bijective (and hence, necessarily birational, since we are
in characteristic 0) morphism of irreducible varieties. In addition
$\ee_{\Omega(n)}$ is nonsingular (see \S 4.1). Thus (4.14) is an isomorphism of varieties. Consequently
(see (4.46)) the map
$$J(n,\ee) \to {\cal O}(M_{\Omega}(n,\ee,\b)),\qquad g\mapsto g|M_{\Omega}(n,\ee,\b)\eqno (4.61)$$ is
an algebra isomorphism. Consequently given any $h\in {\cal O}(M_{\Omega}(n,\ee))$ there exists a
unique $g\in J(n,\ee)$ such that $g|M_{\Omega}(n,\ee,\b) = h|M_{\Omega}(n,\ee,\b)$. But now assume
that $h$ poisson commutes with any function in $J(n,\ee)$. Then $h|M_{\nu}(n,\ee)$ is constant, by
Proposition 4.21, for any $\nu\in \ee_{\Omega(n)}$. But then $h=g$ on $M_{\nu}(n,\ee)$ by (4.25). Thus
$h=g$, by (4.21), and hence $J(n,\ee)$ is maximally Poisson commutative in ${\cal
O}(M_{\Omega}(n,\ee))$. QED\vs 4.8. We recall some definitions, results and notations in Part 1.
Generators $p_{(i)},\,i\in I_{d(n)}$, of the polynomial ring $J(n)$ (see (2.30)) were defined by
(3.20), recalling (2.38). In particular two sets of generators of $J(n)$ were under consideration in
Part 1; namely the $p_{(i)}$ and the $p_i$ (see (2.5) and (2.3)) where $i\in I_{d(n)}$. From the
discussion preceding (2.38) it follows that for $x\in M(n)$, and $m\in I_n$
$$\hbox{span of}\,\,(dp_{(i)})_x,\,\,\,\hbox{for}\,\,i\in I_{d(m)} = \hbox{span
of}\,\,(dp_{i})_x,\,\,\,\hbox{for}\,\,i\in I_{d(m)}\eqno (4.62)$$ and hence $$\hbox{span
of}\,\,(\xi_{p_{(i)}})_x,\,\,\,\hbox{for}\,\,i\in I_{d(m)} = \hbox{span
of}\,\,(\xi_{p_{i}})_x,\,\,\,\hbox{for}\,\,i\in I_{d(m)}\eqno (4.63)$$ But then, by (4.45), for any
$z\in M_{\Omega}(n,\ee)$, $$\hbox{span
of}\,\,(\xi_{\widehat {p_{(i)}}})_z,\,\,\,\hbox{for}\,\,i\in I_{d(m)} = \hbox{span
of}\,\,(\xi_{\widehat {p_{i}}})_z,\,\,\,\hbox{for}\,\,i\in I_{d(m)}\eqno (4.64)$$ An immeditate
consequence of (4.64), when $m=n-1$, and Proposition 4.21 is \vs {\bf Proposition 4.23.} {\it Let
$\nu\in \ee_{\Omega(n)}$ and let $z\in M_{\nu}(n,\ee)$. Then $(\xi_{\widehat {p_{(i)}}})_z,\,i\in
I_{d(n-1)}$, is a basis of $T_z(M_{\nu}(n,\ee))$.}\vs By definition $\a$ (see \S 3.2) is the (complex)
commutative $d(n-1)$-dimensional Lie algebra of vector fields on $M(n)$ spanned by
$\xi_{p_{(i)}},\,i\in I_{d(n-1)}$. By Theorem 3.4 the Lie algebra $\a$ integrates to a (complex)
analytic Lie group $A\cong \Bbb C^{d(n-1)}$ which operates analytically on $M(n)$. Now let $\widehat
{\a}$ be the
$d(n-1)$-dimensional complex commutative (see  (4.44)) Lie algebra of vector fields on
$M_{\Omega}(n,\ee)$ spanned by
$\xi_{\widehat {p_{(i)}}},\,i\in I_{d(n-1)}$. \vs {\bf Remark 4.24.} If $\nu \in \ee_{\Omega(n)}$
then note that, by Proposition 4.23,  $\widehat {\a}|M_{\nu}(n,\ee)$ is a commutative
$d(n-1)$-dimensional Lie algebra of vector fields on $M_{\nu}(n,\ee)$. \vs Let
$\widehat {A}\,(\cong
\Bbb C^{d(n-1)})$ be a simply connected Lie group with Lie algebra $\widehat {\a}$. Let $$A\to \widehat
{A},\,\,a\mapsto
\widehat {a}\eqno (4.65)$$ be the group isomorphism whose differential maps $\xi_{p_{(i)}}$ to
$\xi_{\widehat {p_{(i)}}}$ for all $i\in I_{d(n-1)}$. \vs {\bf Theorem 4.25.} {\it The Lie algebra
$\widehat {\a}$ integrates to an action of $\widehat {A}$ on $M_{\Omega}(n,\ee)$. Furthermore if $a\in
A$ and
$z\in M_{\Omega}(n,\ee)$ then $$\pi_n(\widehat {a}\cdot z) = a\cdot x\eqno (4.66)$$ where $x=
\pi_n(z)$. Moreover
$M_{\nu}(n,\ee)$ is stable under the action of $\widehat {A}$ for any $\nu\in \ee_{\Omega(n)}$. In fact,
for all $\nu$, $\widehat {A}$ operates transitively on $M_{\nu}(n,\ee)$ so that (4.21) is the
decomposition of $M_{\Omega}(n,\ee)$ into $\widehat {A}$ orbits.}\vs {\bf Proof.} Noting Remark 4.24,
Theorem 4.24 is an immediate consequence of (4.45), the
isomorphism (4.54) and Theorem 3.23. QED \vs {\bf Remark 4.26.} Implicit in Theorem 4.25 and its
proof is the fact that if $\nu\in \ee_{\Omega(n)}$ then the Lie algebra, 
$\widehat {\a}|M_{\nu}(n,\ee)$ of vector fields on $M_{\nu}(n,\ee)$ (see Remark 4.24) integrates to
the group action $\widehat {A}|M_{\nu}(n,\ee)$ on $M_{\nu}(n,\ee)$. \vs One now has an
analogue of Theorem 3.5. (Actually it is analogue of a considerably weaker result than Theorem 3.5 in
that
$M_{\Omega}(n,\ee)$ covers the strongly regular set
$M_{\Omega}(n)$ and not all of $M(n)$.) \vs {\bf Theorem 4.27.} {\it Let $f\in J(n,\ee)$ (see (4.46)).
Then the vector field $\xi_f$ integrates to an action of $\Bbb C$ on $M_{\Omega}(n,\ee)$. In fact if
$\nu \in \ee_{\Omega(n)}$ then $\xi_f|M_{\nu}(n,\ee)$ is tangent to $M_{\nu}(n,\ee)$. Indeed
$$\xi_f|M_{\nu}(n,\ee) \in \widehat {\a}|M_{\nu}(n,\ee)\eqno (4.67)$$ so that (see Remark 4.26) the
action of
$\Bbb C$ stabilizes $M_{\nu}(n,\ee)$.}\vs {\bf Proof.} Clearly $\widehat {p_{(i)}}\in J(n,\ee)$, for
$i\in I_{d(n-1)}$, by (4.50). Let $\nu\in \ee_{\Omega(n)}$. Then $\xi_{\widehat
{p_i}}|M_{\nu}(n,\ee),\,\,i\in I_{d(n-1)}$, is a basis of $\widehat {\a}|M_{\nu}(n,\ee)$ by (the
constancy of the $f_i$ in) Theorem 4.22 and Proposition 4.23. But then one has (4.67), also by Theorem
4.22. Theorem 4.27 then follows from Remark 4.26 and Theorem 4.25. QED\vskip 1.5 pc

\centerline{\bf 5. The Emergence of the dual coordinates
$s_j,\,j\in I_{d(n-1)}$}\vskip 1.2pc
5.1. We first wish to be more explicit about the vector fields
$\xi_{r_j},\,j\in I_{d(n-1)}$ on $M_{\Omega}(n,\ee)$. See \S 4.2
and (4.47). Fix $m\in I_n$. We have put $I_{[m]}= I_{d(m)}-I_{d(m-1)}$. For $i\in I_{[m]}$  and $$i= d(m-1) +
k\eqno (5.1)$$ for $k\in I_m$ one has, on $M_{\Omega}(n,\ee)$,
$$\widehat {p_{(i)}} = {1\over m+1 -k}\,\sum_{j\in I_{[m]}} r_j^{m+1 - k}\eqno (5.2) $$ See (2.38), (3.20), \S 4.2
and (4.47). Thus $$d\widehat {p_{(i)}} = \sum_{j\in I_{[m]}} r_j^{m - k}\,dr_j\eqno (5.3)$$ and hence $$
\xi_{\widehat {p_{(i)}}} = \sum_{j\in I_{[m]}} r_j^{m - k}\,\xi_{r_j}\eqno (5.4)$$ (see (1.14)). 

Now let $z\in M_{\Omega}(n,\ee)$ and let $\nu = \kappa_n(z)$ (see (4.9)) so that $z\in M_{\nu}(n,\ee)$. Let $x =
\pi_n(z)$ so that
$x\in M_{\Omega}(n)$. Since the numbers (the eigenvalues of $x_m$) $r_j(z), j\in I_{[m]}$, are distinct the van der
Monde
$m\times m$ matrix $$C_{k\,\ell} = r_j(z)^{m-k} \eqno (5.5)$$ where $$j = d(m-1) +\ell\eqno (5.6)$$ is
invertible.
\vs {\bf Remark 5.1.} In the notation of (4.53) note that $r_j(z) = \nu_{\ell\,m}$. \vs  If $m = n$ then
$\xi_{p_i}$ and $\xi_{p_{(i)}}$ vanish by (2.7) and the argument which implies (2.7). Henceforth assume $m\in
I_{n-1}$. Recalling the definition of the
$m$-dimensional commutative Lie algebra
$\a(m)$ of vector fields on
$M(n)$ (see \S 3.1 and (3.20)) there exists a unique basis $\eta_{j\,\nu},\,j\in I_{[m]},$ of $\a(m)$ 
such that, for $i\in I_{[m]}$, and $k$ related to $i$ by (5.1),
$$\xi_{p_{(i)}} = \sum_{j\in I_{[m]}} r_j^{m - k}(z)\,\eta_{j\,\nu}\eqno (5.7)$$ But then (4.46), (5.4), (5.7) and
the invertibilty of the van der Monde matrix $C_{k\,\ell}$ implies \vs {\bf Proposition 5.2.} {\it Let $\nu \in
\ee_{\Omega(n)}$ and $m\in I_{n-1}$. Then $$(\pi_n)_*(\xi_{r_j}|M_{\nu}(n,\ee)) = \eta_{j\,\nu}|M_c(n)\eqno
(5.8)$$ where $c = \beta(\nu)$ (see (4.22)) for all $j\in I_{[m]}$.}\vs We recall in \S 2.4 that $Z_{x,m}$ is the
commutative (associative) algebra of $M(m)$ generated by $x_m$. Here $x = \pi_n(z)$ so that $x_m$ is regular
semisimple and hence
$dim\,Z_{x,m} = m$. We recall (see \S 3.1) that $G_{x,m}$ is the (algebraic) subgroup of $Gl(n)$
corresponding to $Z_{x,m}$ when $Z_{x,m}$ is regarded as a Lie subalgebra of $M(n)$. We recall also that
$A(m)$ is a simply connected group corresponding to $\a(m)$ and $\a(m)$ integrates to an action of $A(m)$
on $M(n)$ (see Theorem 3.3). Next we recall (see (3.6)) that $\rho_{x,m}$ is the homomorphism of $A(m)$
into $G_{x,m}$ whose differential is given by $$(\rho_{x,m})_*(\xi_{p_{(i)}}) = - (x_m)^{m-k}\eqno (5.9)$$ when
$i\in I_{[m]}$ has the form (5.1). But then applying $(\rho_{x,m})_*$ to (5.7) one has $$-(x_m)^{m-k} =\sum_{j\in
I_{[m]}} r_j^{m - k}(z)\,(\rho_{x,m})_*(\eta_{j\,\nu})\eqno (5.10)$$ On the other hand put $$h_{\nu,m}=
diag(\nu_{1\,m},\ldots, \nu_{m\,m})$$ recalling (4.53), so that in the notation of \S 4.1 $$h_{\nu,m} \in \ee(m)\eqno
(5.11)$$ Note that for the matrix units $e_{\ell\,\ell},\,\ell\in I_m$, one has $$h_{\nu,m}= \sum_{\ell\in
I_m}\,\nu_{\ell\,m}\,e_{\ell\,\ell}\eqno (5.12)$$ Now let $g_z\in Gl(m)$ be such that $$g_z\,h_{\nu,m}g_z^{-1} =
x_m\eqno (5.13)$$ (see Remark 5.1). Using the notation of \S 4.1 note then that $$g_z\,\d(m) g_z^{-1} = Z_{x,m}\eqno
(5.14)$$ \vskip .5pc {\bf Remark 5.3.} Note that $g_z$ is unique in $Gl(m)$ modulo the maximal diagonal torus
$Diag(m)$ (see \S 3.4).\vs  For $j\in I_{[m]}$ let $\varepsilon_{z,j}\in Z_{x,m}$ be the idempotent in
$Z_{x,m}$ defined by putting
$$\varepsilon_{z,j} = g_z\,e_{\ell\,\ell}\, g_z^{-1}\eqno (5.15)$$ where $\ell$ is defined by (5.6). Thus (5.12)
and (5.13) imply that $$x_m = \sum_{\ell\in I_m}\,\nu_{\ell\,m}\,\varepsilon_{z,d(m-1) + \ell}\eqno (5.16)$$ But the
$\varepsilon_{z,d(m-1) + \ell},\,\ell\in I_m$, are orthogonal idempotents by (5.15). Hence, by Remark 5.1,
$$-(x_m)^{m-k} = -\sum_{\ell\in I_m}\,\nu_{\ell\,m}^{m-k}\,\varepsilon_{z,d(m-1) + \ell}\eqno (5.17)$$ \vskip .5pc
{\bf Proposition 5.4.} {\it Let $z\in M_{\Omega}(n,\ee)$ and let $\nu = \kappa(z)$ (see (4.9)). Let $x= \pi_n(z)$
(see (4.8)) and let
$m\in I_n$. Let $\eta_{j,\nu}\in \a(m),\,j\in I_{[m]}$, be the basis of $\a(m)$ defined by (5.7) so that one has
(5.8). Let $(\rho_{x,m})_*:\a(m)\to Z_{x,m}$ be the Lie algebra homomorphism defined as in (3.6) and (3.7). Let
$\varepsilon_{z,j},\,j\in I_{[m]}$, be the orthogonal idempotents in $Z_{x,m}$ defined by (5.15). Then
$$(\rho_{x,m})_*(\eta_{j,\nu}) = - \varepsilon_{z,j}\eqno (5.18)$$ for any $j\in I_{[m]}$.}\vs {\bf Proof.}
 One has $\nu_{\ell\,m}^{m-k} = r_j^{m-k}(z)$ by Remark 5.1. But then (5.18) follows from the equality of the
right hand sides of (5.10) and (5.17), recalling the invertibility of the van der Monde matrix (5.5). QED\vs 5.2.
Retain the notation of Proposition 5.4. For any $\zeta\in \Bbb C^{\times}$ and $i\in I_n$ let $\delta_i(\zeta)\in
Diag(n)$ (see \S 3.4) be the invertible $n\times n$ diagonal matrix such that $\alpha_{j\,j}(\delta_i(\zeta)) =1$
if $j\neq i$ and $\alpha_{i\,i}(\delta_i(\zeta)) = \zeta$, using the notation of (1.2). One notes that if $\ell \in
I_m$ then $\delta_{\ell}(\zeta)\in Diag(m)$. Also, using the relation between $j$ and $\ell$ given by (5.6), 
one has $\gamma_{z\,j}(\zeta)\in G_{x,m}$, by (5.14), where we put
$$\gamma_{z\,j}(\zeta) = g_z\,\delta_{\ell}(\zeta)\,g_z^{-1}\eqno (5.19)$$ Let
$q\in \Bbb C$. One notes that 
$$exp\,\,q\,e_{\ell\,\ell} = \delta_{\ell}(e^q)\eqno (5.20) $$ and hence $$exp\,q\,\varepsilon_{z,j} =
 \gamma_{z\,j}(e^q)\eqno (5.21)$$
Multiplying (5.18) by $q$ and exponentiating (where $exp\,\,q\,\eta_{j,\nu}\in A(m)$) it follows then from
(3.6) that
$$\rho_{x,m}(exp\,\,q\,\eta_{j,\nu}) = \gamma_{z\,j}(e^{-q})\eqno (5.22)$$ \vskip 1pc We can now describe the
flow generated by $\xi_{\r_j}$ (see (4.47) and \S 4.2) on $M_{\Omega}(n,\ee)$ for any $j\in I_{d(n-1)}$ (see Theorem
4.27).\vs {\bf Theorem 5.5.} {\it Let $j\in I_{d(n-1)}$ and let $m\in I_{n-1}$ be such that $j\in I_{[m]} =
I_{d(m)}-I_{d(m-1)}$. Let
$z\in M_{\Omega}(n,\ee)$ (see (4.7)) and let $x\in M_{\Omega}(n),\,\nu\in \ee_{\Omega(n)}$ be such that $x=
\pi_n(z)$ and $\nu =\kappa_n(z)$. See (4.8) and (4.9). Let $q\in \Bbb C$. Then $(\nu,
Ad\,(\gamma_{z\,j}(e^{-q}))(x))\in M_{\nu}(n,\ee)$ (see (4.7) and (4.21)) and $$(exp\,\,q\,\xi_{r_j})\cdot z = (\nu,
Ad\,(\gamma_{z\,j}(e^{-q}))(x))\eqno (5.23)$$ where $\gamma_{z\,j}(e^{-q})\in G_{x,m}$ is defined by
(5.21).}\vs {\bf Proof.} Let $c = \beta(\nu)$ (see (4.22)) so that $c\in \Omega(n)$ and, by Proposition 4.20, the
restriction of $\pi_n$ to $M_{\nu}(n,\ee)$ defines an algebraic isomorphism $M_{\nu}(n,\ee)\to M_c(n)$. But then by
(5.8) one has $$\pi_n((exp\,\,q\,\xi_{r_j})\cdot z) = (exp\,\,q\,\eta_{j,\nu})\cdot x\eqno (5.24)$$ (see Theorems
3.23 and 4.27). But $(exp\,\,q\,\eta_{j,\nu})\in A(m)$ since by definition $\eta_{j,\nu}\in \a(m)$. Hence, by (5.22)
and Theorem 3.3, one has $(exp\,\,q\,\eta_{j,\nu})\cdot x = Ad\,(\gamma_{z\,j}(e^{-q}))(x)$. But then (5.23) follows
from (4.7). QED\vs 5.3. We continue with the notation of \S 5.1 and \S 5.2. \vs {\bf Proposition 5.6.} {\it Let
$j\in I_{d(n-1)}$. Then the isomorphism $exp\,\,q\,\xi_{r_j}$ of $M_{\Omega}(n,\ee)$ reduces to the identity if
$q\in 2\,\pi\,i\,\Bbb Z$.}\vs {\bf Proof.} This is immediate from (5.23) since $\gamma_{z\,j}(e^{-q})$ is the
identity matrix of $M(n)$, by (5.15) and (5.21), if $q\in 2\,\pi\,i\,\Bbb Z$. QED\vs Let $\r$ be the commutative
$d(n-1)$ dimensional Lie algebra of vector fields on $M_{\Omega}(n,\ee)$ with basis $\xi_{r_j},\,j\in I_{d(n-1)}$
(see (4.57) and (4.59)). Let $j\in I_{d(n-1)}$ and let $A_{\r,j}$ be a 1-dimensional complex torus with a global
coordinate $\zeta_j$ defining an algebraic group isomorphism $$\zeta_j:A_{\r,j}\to \Bbb C^{\times}\eqno (5.25)$$ By
Proposition 5.6 (and abuse of notation) we can regard $$ \Bbb C\,\xi_{r_j} = Lie\,A_{\r,j}\eqno (5.26)$$ and
simultaneously have
$A_{\r,j}$ operate on $M_{\Omega}(n,\ee)$, as an integration of the vector field $\xi_{r_j}$, in such a fashion that
if
$b\in A_{\r,j}$ then
$b = exp\,q\,\xi_{r_j}$ if
$$ \zeta_j(b) = e^{q}\eqno (5.27)$$ It is very easy to prove that the action of $A_{\r,j}$ is analytic but what is
much more important for us is to prove that this action is that of an algebraic group, operating algebraically on 
an affine algebraic variety. \vs {\bf Theorem 5.7.} {\it Let $j\in I_{d(n-1)}$. Then the map $$A_{\r,j}\times
M_{\Omega}(n,\ee)\to M_{\Omega}(n,\ee),\qquad (b,z)\mapsto b\cdot z \eqno (5.28)$$ is a (algebraic) morphism.}\vs
{\bf Proof.} Let $b\in A_{\r,j}$ and let $z\in M_{\Omega}(n,\ee)$. Let $\zeta = \zeta_j(b)$ so that $\zeta\in
\Bbb C^{\times}$. Recalling (4.7) write $z = (\nu,x)$ where $x\in M_{\Omega}(n)$ and $\nu\in \ee_{\Omega(n)}$.
Let $m\in I_{n-1}$ be such that $j\in I_{[m]} = I_{d(m)}- I_{d(m-1)}$. Then, by (5.23) and (5.27), $$b\cdot z =
 (\nu,Ad\,(\gamma_{z,j}(\zeta^{-1}))(x))\eqno (5.29)$$ Since $\nu = \kappa_n(z)$ and $x=\pi_n(z)$ (see (4.8) and
(4.9)) and both (4.8) and (4.9) are morphisms it suffices by (5.29) to prove that the map
 $$A_{\r,j}\times M_{\Omega}(n,\ee)\to G_{x,m},\qquad (b,z) \mapsto \gamma_{z,j}(\zeta^{-1})\eqno (5.30)$$ is a
morphism. Now $\ell\in I_m$ in (5.19) is defined so that $j= d(m-1) + \ell$. For any $g\in Gl(m)$ let $[g]\in
Gl(m)/Diag(m)$ be the left coset defined by $g$ (see \S 3.4). Of course $Gl(m)/Diag(m)$ is an affine algebraic
homogeneous space. Recalling (5.19) to prove (5.30) is a morphism it clearly suffices to show that $$
A_{\r,j}\times M_{\Omega}(n,\ee) \to Gl(m)/Diag(m),\qquad (b,z) \mapsto [g_z]\eqno (5.31)$$ is a
morphism.

Let $E(m)$ be the set of all regular semisimple elements in $M(m)$ so that $E(m)$ has the structure of a Zariski open
(and hence nonsingular) affine subvariety of $M(n)$. Then, recalling (4.1), the map $$Gl(m)\times \ee(m)\to
E(m),\qquad (g,\mu)\mapsto g\,\mu\,g^{-1}\eqno (5.32)$$ is a surjective morphism. But now if $g\in Gl(m)$ and $\mu
\in \ee(m)$ then $[g]\cdot
\mu\in E(m)$ is well defined by putting $[g]\cdot \mu = g\,\mu\,g^{-1}$. Clearly $$(Gl(m)/Diag(m)) \times \ee(m)\to
E(m),\qquad ([g],\mu)\mapsto [g]\cdot \mu \eqno (5.33)$$ is then also a surjective morphism. Now let 
$$E(m,\ee) = \{(\mu,y)\in  \ee(m)\times E(m)\mid \mu\,\,\hbox{is $Gl(m)$-conjugate to $y$}\}\eqno (5.34)$$ so 
that $E(m,\ee)$ is a Zariski closed subset of $\ee(m)\times E(m)$. The argument establishing the dimension and
nonsingularity in Theorem 4.9 (especially using the independence of the differentials $dp_i,\,i\in I_{[m]}$, at
all points in $E(m)$) can obviously be modified to apply here and prove that
$$E(m,\ee)\,\,\hbox{is a nonsingular
$m^2$-dimensional Zariski closed subset of $\ee(m)\times E(m)$}\eqno (5.35)$$ But now (5.33) may be augmented to
define the map $$(Gl(m)/Diag(m)) \times \ee(m)\to E(m,\ee),\qquad ([g],\mu)\mapsto (\mu,[g]\cdot \mu )\eqno (5.36)$$
But (5.33) readily implies that (5.36) is a surjectve morphism so that, for one thing, $E(m,\ee)$ is irreducible.
Hence $E(m,\ee)$ is a nonsingular variety. But (5.36) is obviously bijective and hence birational. But then 
(5.36) is an algebraic isomorphism. Let $$E(m,\ee) \to (Gl(m)/Diag(m)) \times \ee(m)\eqno (5.37)$$ be
the inverse isomorphism. Projecting on the first factor defines a morphism $\sigma:E(m,\ee)\to Gl(m)/Diag(m)$
where, for $g\in Gl(m)$ and $\mu \in \ee(m)$, $$\sigma((\mu,[g]\cdot \mu)) = [g]\eqno (5.38)$$ But now since
(4.8) and (4.9) are morphisms it follows that $\tau:A_{\r,j}\times M_{\Omega}(n,\ee)\to E(m,\ee)$ is a morphism
where, using the noataion of (4.3), $\tau ((b,z)) = (\nu(m),x_m)$. But $$\sigma\circ \tau ((b,z)) = [g_z]\eqno
(5.39)$$ by (5.13), since clearly $h_{\nu,m} = \nu(m)$. See (4.3) and (5.12). This proves that (5.31) is a morphism.
QED\vs 5.4. Let $m\in I_{n-1}$ and let $\r(m)$ be the span of the vector fields $\xi_{r_i},\,i\in
I_{[m]}$, and so that, as defined in \S 5.3, $$\r = \r(1) \oplus \cdots \oplus \r(m-1)\eqno (5.40)$$  By (4.57) and
(4.59),
$\r(m)$ is a commutative Lie algebra of dimension $m$ and as we have already noted $\r$ is a commutative Lie
algebra of dimension
$d(n-1)$. Let (see (5.26))
$$\eqalign{A_{\r}(m) &= A_{\r,d(m-1)+1}\times \cdots \times A_{\r,d(m-1) + m}\cr A_{\r} &= A_{\r}(1)\times \cdots
\times A_{\r}(n-1)\cr}\eqno (5.41)$$ so that as algebraic groups $$\eqalign{A_{\r}(m) &\cong (\Bbb C^{\times})^{m}\cr
A_{\r}&
\cong (\Bbb C^{\times})^{d(n-1)}}\eqno (5.42)$$ In addition, by (5.26), $$\eqalign{\r(m)&= Lie\,A_{\r}(m)\cr \r&=
Lie\,A_{\r}\cr}\eqno (5.43)$$ As an immediate consequence of Theorem 5.7 and commutativity one has \vs {\bf Theorem
5.8.} {\it Let $m\in I_{n-1}$. Then the Lie algebras $\r(m)$ and $\r$, respectively, integrate to an algebraic
action of $A_{\r}(m)$ and $A_{\r}$ on $M_{\Omega}(n,\ee)$.}\vs The following result is a refinement of Theorem
4.25. We are now dealing with the ``eigenvalue" vector fields $\xi_{r_i}$ themselves on $M_{\Omega}(n,\ee)$ rather
than the more crude ``eigenvalue symmetric function" vector fields $\xi_{\widehat {p_{(i)}}}$.\vs {\bf Theorem 5.9.}
{\it Let
$\nu \in \ee_{\Omega(n)}$. Then $M_{\nu}(n,\ee)$ is stable under the algebraic group $A_{\r}$. Furthermore $A_{\r}$
operates simply and transitively on $M_{\nu}(n,\ee)$. In particular the disjoint union (4.21) is the $A_{\r}$-orbit
decomposition of $M_{\Omega}(n,\ee)$.}\vs {\bf Proof.} By Theorem 4.22 and Theorem 4.27 one has $$\r|M_{\nu}(n,\ee)=
\widehat {\a}|M_{\nu}(n,\ee)\eqno (5.44)$$ Hence $$A_{\r}|M_{\nu}(n,\ee) = \widehat {A}|M_{\nu}(n,\ee)\eqno
(5.45)$$ But then
$M_{\nu}(n,\ee)$ is stable under
$A_{\r}$ and
$A_{\r}$ operates transitively on $M_{\nu}(n,\ee)$ by Theorem 4.25. The only question concerns the simplicity of
this action.

Let $b\in A_{\r}$. By definition there exists $q_j\in \Bbb C,\,j\in I_{d(n-1)}$, such that 
$$b = exp\,q_1\,\xi_{r_1}\cdots exp\,q_{d(n-1)}\,\xi_{r_{d(n-1)}}\eqno (5.46)$$ But then if $c = \beta(\nu)$ (see
(4.22)) it follows from Proposition 5.2 that there exists $\eta_{j\,\nu}\in \a,\,j\in I_{d(n-1)}$, such that if $a\in
A$ is defined by putting $$a = exp\,q_1\,\eta_{1,\nu}\cdots exp\,q_{d(n-1)}\,\eta_{d(n-1),\nu}\eqno (5.47)$$ then
$$b|M_{\nu}(n,\ee) =
\widehat {a}|M_{\nu}(n,\ee)\eqno (5.48)$$ recalling Theorem 4.25. In addition Proposition 5.2 implies that $a =
a(1)\cdots a(n-1)$ where, for $m\in I_{n-1}$, $a(m) \in A(m)$ is given by $$a(m) = \prod_{j\in I_{[m]}}
exp\,q_{j}\,\eta_{j,\nu}\eqno (5.49)$$  Now assume that
$b|M_{\nu}(n,\ee)$ has a fixed point. But then by the commutativity of $A_{\r}$ and the transitivity of
$A_{\r}$ on $M_{\nu}(n,\ee)$ it follows that $b|M_{\nu}(n,\ee)$ reduces to the identity. We must prove $$q_j\in
2\,\pi\,i\,\Bbb Z,\,\,\forall j\in I_{d(n-1)}\eqno (5.50)$$ by (5.27) and (5.41). But now
$a|M_{c}(n)$ reduces to the identity by (4.66) and (5.48). See Proposition 4.20. But then $a\in D_c$ (see
\S 3.64). Hence $a(m)\in D_c(m)$ for any $m\in I_{n-1}$ by Theorem 3.28. But Theorem 3.28 also asserts that
if $z\in M_{\nu}(n,\ee)$ and $x = \pi_n(z)$ then also $a(m)\in Ker\,\rho_{x,m}$. But by (5.15) and
(5.18) (see also (5.19) and (5.22)) this implies (5.50) since $m\in I_{n-1}$ is arbitrary. QED\vs 5.5. In the
introduction, \S 0, we defined $\b_e\s M(n)$. In common parlance (for some people) $\b_e$ is the space of all
$n\times n$
 Hessenberg
matrices. In \S 4.3 we defined $\b_{e,\Omega(n)}$ to be the intersection $M_{\Omega}(n)\cap \b_e$ so that 
$\b_{e,\Omega(n)}$ is a Zariski open subvariety of $\b_e$ (and hence a nonsingular variety) and a closed subvariety
of $M_{\Omega}(n)$. We also defined $M_{\Omega}(n,\ee,\b) = \pi_{n}^{-1}(\b_{e,\Omega(n)})$ (see (4.8) and \S4.3) so
that 
$M_{\Omega}(n,\ee,\b)$ is a Zariski closed subset of $M_{\Omega}(n,\ee)$. Sharpening Theorem 4.4 we
shall need
\vs {\bf  Theorem 5.10.} {\it The restriction (see (4.9)) $$\kappa_n:M_{\Omega}(n,\ee,\b)\to \ee_{\Omega(n)}\eqno
(5.51)$$ is an algebraic isomorphism so that $M_{\Omega}(n,\ee,\b)$ is a closed nonsingular subvariety of
$M_{\Omega}(n,\ee)$.}\vs {\bf Proof.} Clearly (5.51) is a morphism since it the restriction of (4.9) to a Zariski 
closed subset of $M_{\Omega}(n,\ee)$. On the other hand it is bijective by Theorem 4.4. However the inverse of (5.51)
is a morphism. See (4.20). QED\vs An easy consequence of Theorem 5.10 is \vs {\bf Theorem 5.11.} {\it The image of
the map
$$M_{\Omega}(n,\ee,\b) \to \Bbb C^{d(n)},\qquad y\mapsto (r_1(y),\ldots,r_{d(n)}(y))\eqno (5.52)$$ is a Zariski
open set in $\Bbb C^{d(n)}$ of the form $(\Bbb C^{d(n)})_q$ where $q$ is a nonzero polynomial on $\Bbb C^{d(n)}$ and
(5.52) is an algebraic isomorphism of
$M_{\Omega}(n,\ee,\b)$ with $(\Bbb C^{d(n)})_q$.} \vs {\bf Proof.} Recalling the definition of $\ee_{\Omega(n)}$ in
\S 4.1, Theorem 5.11 follows immediately from Theorem 5.10 and the definition of $r_i$ in (4.47) and
$\rho_{k\,m}$ in \S 4.2. QED \vs On the other hand 
we establish the following product structure for
$M_{\Omega}(n,\ee)$.
\vs {\bf Theorem 5.12.} {\it The map $$A_{\r}\times M_{\Omega}(n,\ee,\b)
\to M_{\Omega}(n,\ee),\qquad (b,y)\mapsto b\cdot y \eqno (5.53)$$ is an algebraic isomorphism.} \vs {\bf Proof.} The
map (5.53) is bijective by (4.21), (4.25) and Theorem 5.9. But then (5.53) is a bijective morphism of nonsingular
algebraic varieties by Theorem 5.8 and Theorem 5.10. Hence (5.53) is an algebraic isomorphism. QED \vs By
definition (see (5.41)) any element $b\in A_{\r}$ can be uniquely written $$b = (b_1,\ldots,b_{d(n-1)})\eqno
(5.54)$$ where $b_j\in A_{\r,j}$. We now extend the domain of the function $\zeta_j$ on $A_{\r,j}$ to all of
$A_{\r}$ so that, in the notation of (5.54), $$\zeta_j(b) =\zeta_j(b_j)\eqno (5.55)$$ Thus $\zeta_j\in {\cal
O}(A_r)$ and the map $$A_{\r}\to (\Bbb C^{\times})^{d(n-1)},\qquad b\mapsto
(\zeta_1(b),\ldots,\zeta_{d(n-1)}(b))\eqno (5.56)$$ is an isomorphism of algebraic groups. For $i\in d(n-1)$ let
$\lambda_{r_i}$ be the left invariant vector field on $A_{\r}$ whose value at the identity of $A_{\r}$ corresponds
to (abuse of notation) $\xi_{r_i}$. Thus, by (5.27), in the coordinates $\zeta_j$ of $A_{\r}$, one has
$$\lambda_{r_i} = \zeta_i\,{\partial\over \partial\,\zeta_i}$$ and hence $$\lambda_{r_i}\,\zeta_j=
\delta_{i\,j}\zeta_i\eqno (5.57)$$\vskip .5pc But now, by Theorem 5.12, every $z\in M_{\Omega}(n,\ee)$ can be
uniquely written $z = b\cdot y$ where $b\in A_{\r}$ and $y\in M_{\Omega}(n,\ee,\b)$. Hence, by Theorem 5.12, one has
a well--defined function $s_{j}\in {\cal O}(M_{\Omega}(n,\ee)),\,j\in I_{d(n-1)}$, where $$s_j(z) =
\zeta_j(b^{-1})\eqno(5.58)$$ Furthermore Theorems 5.8, 5.9, 5.12 and (5.55) also clearly imply \vs {\bf Theorem 5.13.}
{\it Let $\nu\in 
\ee_{\Omega(n)}$. Then the map $$M_{\nu}(n,\ee) \to (\Bbb C^{\times})^{d(n-1)},\qquad z\mapsto
(s_1(z),\ldots,s_{d(n-1)}(z))$$ is an algebraic isomorphism.}\vs The action of $A_{\r}$ on
$M_{\Omega}(n,\ee)$ of course introduces, contragrediently, an action of $A_{\r}$ on ${\cal O}(M_{\Omega}(n,\ee))$ so
that if $b\in A_{\r},\,\, f\in {\cal O}(M_{\Omega}(n,\ee)$ and $z\in M_{\Omega}(n,\ee)$ then $(b\cdot f)(z) =
f(b^{-1}\cdot z)$. It is immediate from (5.58) that $$b\cdot s_j = \zeta_j(b)\,s_j,\,\,j\in I_{d(n-1)}\eqno
(5.59)$$ But then with regard to Poisson bracket on
$M_{\Omega}(n,\ee)$, at this stage, we can say, for $i\in I_{d(n)}$ and $j\in I_{d(n-1)}$, $$[r_i,s_j] =
\delta_{i\,j} s_j\eqno (5.60)$$ since, by differentiating (5.59) and applying (5.57) clearly,
$$\xi_{r_i}s_j = \delta_{i\,j} s_j\eqno (5.61)$$ recalling (4.59). 

Combining Theorems 5.11, 5.12, 5.13 and (5.61) one sees that the $s_i,r_j,\,i\in I_{d(n-1)},\,j\in I_{d(n)}$, form a
system of uniformizing parameters on $M_{\Omega}(n,\ee)$. \vs {\bf Theorem 5.14.} {\it For any $z\in
M_{\Omega}(n,\ee)$ the $n^2$ differentials $$(dr_i)_z,\,(ds_j)_z,\,i\in I_{d(n)},\,j\in I_{d(n-1)},\,\,\hbox{are a
basis of the cotangent space $T^*_z(M_{\Omega}(n,\ee))$}\eqno (5.62)$$\vskip .5pc Furthermore the image of the map
$$M_{\Omega}(n,\ee)\to \Bbb C^{n^2},\qquad z \mapsto (r_1(z),\ldots,r_{d(n)}(z),s_1(z),\ldots,s_{d(n-1)}(z))
\eqno (5.63)$$ is a Zariski open set $Y$ in $\Bbb C^{n^2}$ and (5.63) is an algebraic isomorphism of $M_{\Omega}(n,\ee)$
with $Y$.} \vs {\bf Proof.} Let $z\in M_{\Omega}(n,\ee)$. The $(dr_i)_z,\,i\in I_{d(n)},$ are linearly
independent by Proposition 4.18 (double prime statement). On the other hand $(ds_j)_z,\,i\in I_{d(n)},\,j\in
I_{d(n-1)}$, are clearly linearly independent by (5.61). But (5.61) together with (4.57) implies that $(dr_i)_z$
are independent of the $(ds_j)_z$. By dimension this proves (5.62). The image of (5.63) is the just the product
of the image of (5.52) and (5.59) by Theorem 5.12. Here we are using the constancy of the $r_i$ on $M_{\nu}(n,\ee)$
(see (4.53)). But then the image $Y$ of (5.63) is Zariski open in $\Bbb C^{n^2}$ by Theorems 5.11 and  5.13. But the
morphism (5.63) is bijective by Theorems 5.11, 5.12 and 5.13. Since both $Y$ and $M_{\Omega}(n,\ee)$ are nonsingular
varieties it follows (a version of Zariski's Main Theorem) that (5.63) is an algebraic isomorphism. QED\vs 5.6. It
is our main objective now to prove the $s_i,\,i\in I_{d(n-1)}$, Poisson commute among themselves. Localizing this problem, 
letting $z\in M_{\Omega}(n,\ee)$ and $i,j\in I_{d(n-1)}$ it is enough to show that $$[s_i,s_j](z) = 0 
\eqno (5.64)$$ Let $\nu_o = \kappa_n(z)$ so that $z\in M_{\nu_o}(n,\ee)$. By (4.25) there exists a unique element
$z_o$ in $M_{\Omega}(n,\ee,\b)\cap M_{\nu_o}(n,\ee)$. Let $x_o= \pi_n(z_o)$ (see (4.8)). Recalling the definition of
$M_{\Omega}(n,\ee,\b)$ and $\b_{e,\Omega(n)}$ in \S 4.3 it follows that $x_o\in \b_{e,\Omega(n)}$ and, by
Proposition 4.14, $\pi_n$ defines $M_{\Omega}(n,\ee,\b)$ as an analytic covering space of $\b_{e,\Omega(n)}$ (see
also Theorem 5.10). Thus there exists an open connected (in the Euclidean sense) neighborhood $V_{x_o}$ of $x_o$ in 
$\b_{e,\Omega(n)}$ such that $V_{x_o}$ is evenly covered by $\pi_n$. Let $V_{z_o}$ be the connected component of 
$\pi_n^{-1}(V_{x_o})$ which contains $z_o$. Thus $V_{z_o}$ is an open connected neighborhood of $z_o$ in
$M_{\Omega}(n,\ee,\b)$ and $$\pi_n:V_{z_o}\to V_{x_o}\eqno (5.65)$$ is an analytic isomorphism. 

Now, recalling the analytic isomorphism (4.13), let $c_o= \Phi_n(x_o)$ (see (2.8)) and let $V_{c_o}$ be 
the open connected neighborhood of $c_o$ in $\Omega(n)$ defined by putting $V_{c_o} = \Phi_n(V_{x_o})$. Finally
let $V_{\nu_o} = \kappa_n(V_{z_o})$ so that (see (5.51)) $V_{\nu_o}$ is open connected neighborhood of $\nu_o $ in
$\ee_{\Omega(n)}$. \vs {\bf Lemma 5.15.} {\it The map $$\beta \circ \kappa_n: V_{z_o}\to V_{c_o}\eqno (5.66)$$ is
an analytic isomorphism (see (4.15) for the definition of $\beta$).}\vs {\bf Proof.} One readily notes that $\beta
\circ \kappa_n$ restricted to $V_{z_o}$ is the same as $\Phi_n\circ \pi_n$ restricted to $V_{z_o}$ (a commutative
diagram). But (5.65) is an analytic isomorphism and (4.13) is an analytic isomorphism. QED\vs Let $x =
\pi_n(z)$. Let $W_z = A_{\r}\cdot V_{z_o}$ and let $W_x = \pi_n(W_z)$. \vs {\bf Proposition 5.16.} {\it $W_z$ is an
open connected neighborhood of $z$ in $M_{\Omega}(n,\ee)$. Furthermore $$W_z = \sqcup_{\nu\in
V_{\nu_o}}\,M_{\nu}(n,\ee) \eqno (5.67)$$ In addition $W_x$ is a open connected neighborhood of $x$ in
 $M_{\Omega}(n)$ and $$W_x = \sqcup_{c\in V_{c_o}}M_c(n)\eqno (5.68)$$ (see \S2.2). Finally $$\pi_n:W_z\to W_x\eqno
(5.69)$$ is an analytic isomorphism.}\vs {\bf Proof.} Since $V_{z_o}$ is open and connected in
$M_{\Omega}(n,\ee,\b)$ it follows from Theorem 5.12 that $W_z$ is open and connected in $M_{\Omega}(n,\ee)$.
Furthermore (5.67) follows from Theorems 5.9 and 5.12. Also $z\in W_z$ since $\nu_o\in V_{\nu_o}$. Since a covering
map is an open map it follows that $W_x$ is an open connected neighborhood of $x$ in $M_{\Omega}(n)$. Also (5.68)
follows from (5.67) and Proposition 4.20. Since (5.69) is both open and continuous, to prove that it is an analytic
isomorphism it suffices to prove that it is bijective. For this of course one must see that it is injective. But
this clearly follows from Proposition 4.20, Lemma 5.15 and the bijectivity of (5.51). QED \vs Proposition 5.16
enables us to carry holomorphic functions on $W_z$ to $W_x$. Using the notation of \S 1.2, for $i\in I_{d(n)},\,j\in
I_{d(n-1)}$, let $r_i',\,s_j'\in {\cal H}(W_x)$ be defined so that on $W_z$, $r_i'\circ \pi_n = r_i$ and $s_j'\circ
\pi_n= s_j$. Recalling the definition of the Poisson structure on $M_{\Omega}(n,\ee)$ (see \S 4.6) it follows that
(5.69) is an isomorphism of Poisson manifolds. Hence to prove (5.64) it suffices to prove $$[s_i',s_j'](x) = 0\eqno
(5.70)$$\vskip.5pc Obviously, for $i\in I_{d(n-1)}$, $$(\pi_n)_*(\xi_{r_i}|W_z )= \xi_{r_i'}\eqno (5.71)$$ If 
$\r'$ is the commutative Lie algebra span of $\xi_{r_i'},\,i \in I_{d(n-1)}$, then, recalling 
Proposition 4.20, Theorem 5.9, (5.67) and (5.68), $\r'$ integrates to a group $A_{\r'}$, which operates on $W_x$,
and which admits an isomorphism $$A_{\r}\to A_{\r'},\qquad b \mapsto b'\eqno (5.72)$$ such that for any
$w\in W_z$, and $b\in A_{\r}$ $$\pi_n(b\cdot w) = b'\cdot \pi_n(w)\eqno (5.73)$$ In addition Theorems 5.9 and
5.12 imply
\vs {\bf Proposition 5.17.} {\it (5.68) is the orbit decomposition of $A_{\r'}$ on $W_x$ and $A_{\r'}$ operates
simply (as well as transitively) on $M_c(n)$ for every $c\in V_{c_o}$. Furthermore (see (5.65)) $$W_x = \sqcup_{b\in
A_r}\,\, b'\cdot  V_{x_o}\eqno (5.74)$$}.\vs 5.7. We retain the notation of the previous section and recall some
notation from
\S 1.2. If
$O\s M(n)$ is an adjoint orbit of $Gl(n)$ and $y\in O$ then $O_y = O$. If $O$ is an orbit of regular semisimple
elements then
$$dim\,O = 2\,d(n-1)\eqno (5.75)$$ and if $q_k,\,k=1,\ldots,n$, is the constant value that the
$Gl(n)$-invariant (see \S 2.1) $p_{d(n-1)+k}$ takes on $O$ then $O$ is clearly determined by those values. That
is, $$O =
\{y\in M(n)\mid p_{d(n-1)+k}(y) = q_k\}\eqno (5.76)$$ It follows then from (2.9) that if (see \S
2.2) $$\Bbb C^{d(n)}(O)= \{c\in \Bbb C^{d(n)}\mid c_{d(n-1) + k} = q_k, \,k\in I_n\} $$ one has $$O = \sqcup_{c\in 
\Bbb C^{d(n)}(O)}\, M_c(n)\eqno (5.77)$$ Let ${\cal S}$ be the set of all $Gl(n)$-adjoint orbits of regular
semisimple elements. Since any $y\in M_{\Omega}(n)$ (see (2.53)) is regular semisimple one has $O_y\in {\cal S}$
for any $y\in W_x$. Let $${\cal S}(W_x) = \{O\in {\cal S}\mid W_x\cap O \neq \emptyset\}$$ and for any $O\in {\cal
S}(W_x)$ let $V_{c_o}(O) = \Bbb C^{d(n)}(O)\cap V_{c_o}$ so that $$O\cap W_x = \sqcup_{c\in V_{c_o}(O)} 
M_c(n)\eqno
(5.78)$$ by (5.68) and (5.77).

Now, recall from \S 1.2, any adjoint orbit has the structure of a symplectic manifold. If $O\in {\cal S}$
then as one knows $O$ is closed in $M(n)$. If $O\in {\cal S}(W_x)$ then $O\cap W_x$ is open in $O$ and hence
$O\cap W_x$ is a symplectic manifold. \vs {\bf Proposition 5.18.} {\it Let $O\in {\cal S}(W_x)$. Then the group
$A_{\r'}$ (see Proposition 5.17) stabilizes $O\cap W_x$ and operates as a group of symplectomorphisms on 
$O\cap W_x$.} \vs {\bf Proof.} For any $i\in I_{d(n-1)}$ the vector field $\xi_{r_i'}|O\cap W_x$ is tangent to
$O\cap W_x$ and is a Hamiltonian vector field on $O\cap W_x$ by Proposition 1.3 and especially (1.19). Thus
$\r'|O\cap W_x$ is a Lie algebra of Hamiltonian vector fields on $O\cap W_x$. But, by (5.78) and Proposition 5.17,
$O\cap W_x$ is stabilized by the integrated group $A_{\r'}$. Hence $A_{\r'}$ operates as a group of symplectomorphisms
of $O\cap W_x$. QED\vs Assume $X$ is a submanifold of $W_x$. Let ${\cal S}(X) = \{O\in {\cal S}\mid O\cap X\neq 
\emptyset\}$ so that ${\cal S}(X)\s {\cal S}(W_x)$. We will say $X$ is Lagrangian in $W_x$ if $dim\,X = d(n)$ and
$O\cap X$ is a Lagrangian submanifold of the symplectic manfold $O$, for any $O\in {\cal O}(X)$. \vs {\bf Proposition
5.19.} {\it Assume $X$ is a Lagrangian submanifold of $W_x$. Let $b\in A_{\r}$. Then ${\cal S}(X) = {\cal S}(b'\cdot 
X)$ and $b'\cdot X$ is again a Lagrangian submanifold of $W_x$. }\vs {\bf Proof.} Obviously $dim\,\,\,b'\cdot X =
d(n)$. Let
$O\in {\cal S}(X)$. Then of course $O\cap X = (O\cap W_x)\cap X$ and, by definition of being Lagrangian, obviously
$(O\cap W_x)\cap X$ is Lagrangian in
$O\cap W_x$. But $b'$ operates as a symplectomorphism of $O\cap W_x$ by Proposition 5.18. Thus $b'\cdot ((O\cap
W_x)\cap X) = (O\cap W_x)\cap b'\cdot X$ is again Lagrangian in $O\cap W_x$. That is, $O\in {\cal S}(b'\cdot X)$
and $b'\cdot (O\cap X) = O\cap b'\cdot X$ is Lagrangian in $O$. By using $b^{-1}$ one readily reverses the argument to
show that if $O\in {\cal S}(b'\cdot X)$ then $O\in {\cal S}(X)$ and $O\cap b'\cdot X$ is Lagrangian in $O$. QED \vs
The following result will be seen to be the key point in proving (5.70) and consequently (5.64). \vs {\bf Theorem
5.20.}  {\it Retain
the notation of (5.65) (or (5.74)). Then $V_{x_o}$ is a Lagrangian submanifold of $W_x$.} \vs {\bf Proof.} We will
use results in [K2]. These are stated for complex semisimple Lie groups but their extension to the reductive
group $Gl(n)$ is immediate and we will apply the results for that case. Let $N\s Gl(n)$ be the maximal unipotent
subgroup where $Lie\,\,N$ is the Lie subalgebra of all strictly upper triangular matrices. Retaining notation in [KW]
we have put $\u = Lie\,N$. See (3.51)). Let $\ss$ be defined by (1.1.5) in [K2] so that if $\ss_e = -e +\ss$, using
the notation of \S 2.2, then $\ss_e\s \b_e$ (see \S4.3) and, as asserted by Theorem 1.1 in [K2], (a), $\ss_e$ is a
cross-section for the adjoint action of $Gl(n)$ on the set of all regular elements in $M(n)$. On the other hand (b),
Theorem 1.2 in [K2] asserts that
$\b_e$ is stable under
$Ad\,N$ and the map $$N\times \ss_e\to \b_e,\qquad (u,w)\mapsto Ad\,u (w)\eqno (5.79)$$ is an isomorphism of affine
varieties.

To prove that $V_{x_o}$ is Lagrangian in $W_x$ we first observe that $dim\,V_{x_o} = d(n)$. This is clear since
$V_{x_o}$ is open in $\b_e$ (see \S 5.6). Now let $O\in {\cal S}(V_{x_o})$. It remains to show that $O\cap V_{x_o}$
is Lagrangian in $O$. But now, by (a) above, $O\cap \ss_e$ consists of a single point $y_o$ and hence by (b) one must
have $$O\cap \b_e =Ad\,N(y_o)\eqno (5.80)$$ Since $O$ is closed in $M(n)$ this implies that $Ad\,N(y_o)$ is 
closed in $\b_e$ and one has $$O\cap V_{x_o}= Ad\,N(y_o)\cap V_{x_o}\eqno (5.81)$$ since
$V_{x_o}\s \b_e$. But $V_{x_o}$ is open in $\b_e$ and hence (5.81) implies that $O\cap V_{x_o}$ is open in 
$Ad\,N(y_o)$. Consequently to prove that $O\cap V_{x_o}$ is Lagrangian in $O$ it suffices to prove that $Ad\,N(y_o)$
is Lagrangian in $O$. But $dim\,\,Ad\,N(y_o) = dim\,N$ by (5.79), and $dim\,N = d(n-1)$ which is half the dimension
of 
$O$. On the other hand if $y\in Ad\,N(y_o)$ and $u,v\in \u$ (see (3.51)) we must show that $\omega_y(\eta^u,\eta^v)=0$
(see 
\S 1.2 and more specifically (1.10)). But $\omega_y(\eta^u,\eta^v)= B(y,[u,v])$ by (1.10). But $B(y,[u,v])= 0$ since
$y\in \b_e$ and one notes that $\b_e$ is $B$-orthogonal to $[\u,\u]$. QED \vs By Theorem 5.12 one has the following 
disjoint union $$M_{\Omega}(n,\ee) = \sqcup_{b\in A_{\r}}\,b\cdot M_{\Omega}(n,\ee,\b)\eqno (5.82)$$ One the other
hand, we note that the components in (5.82) are level sets of the functions $s_j,\,j\in I_{d(n-1)}$. Indeed for any
$\tau\in (\Bbb C^{\times})^{d(n-1)}$ let $\tau_i\in \Bbb C^{\times},\,i\in I_{d(n-1)}$, be defined so that $$\tau =
(\tau_1,\ldots,\tau_{d(n-1)})\eqno (5.83)$$ and for $\tau \in (\Bbb C^{\times})^{d(n-1)}$, let
$$M_{\Omega}(n,\ee,\tau) = \{y \in M_{\Omega}(n,\ee)\mid s_i(y) = \tau_i,\,\forall i\in I_{d(n-1)}\}\eqno (5.84)$$ so
that $$M_{\Omega}(n,\ee) = \sqcup_{\tau\in (\Bbb C^{\times})^{d(n-1)}}\,M_{\Omega}(n,\ee,\tau)\eqno (5.85)$$ The
following proposition is an immediate consequence of (5.58). \vs {\bf Proposition 5.21.} {\it The partitions (5.74)
and (5.82) of $M_{\Omega}(n,\ee)$ are the same. That is, for any $b\in A_{\r}$, $$b\cdot M_{\Omega}(n,\ee,\b) =
M_{\Omega}(n,\ee,\tau)\eqno (5.86)$$ where for $i\in I_{d(n-1)},\,\,\tau_i= \zeta_i(b^{-1})$.}\vs {\bf Remark 5.22.} Since
$V_{z_o}$ is open in $M_{\Omega}(n,\ee,\tau)$ (see \S 5.6) note from the definition $W_z$ in \S 5.6 one has the
disjoint union $$W_z = \sqcup_{b\in A_{\r}}\,b\cdot V_{z_o}\eqno (5.87)$$ and hence by (5.86), $$ b\cdot
V_{z_o}= \{w\in W_{z}\mid s_i(w)= \zeta_i(b^{-1}),\,\forall i\in I_{d(n-1)}\}\eqno (5.88)$$\vskip .5pc We can now prove
one of the main theorems of the paper. \vs {\bf Theorem 5.23.} {\it The uniformizing parameters $r_i,s_j,\,i\in
I_{d(n)},\,j\in I_{d(n-1)}$ of $M_{\Omega}(n,\ee)$ (see Theorem 5.14) satisfy the following Poisson commutation
relations, $$ \eqalign{(1)&\,\,[r_i,r_j]= 0,\,i,j\in I_{d(n)}\cr (2)&\,\,[r_i,s_j] = \delta_{i\,j}\,s_j,\,i\in
I_{d(n)},\,j\in I_{d(n-1)}\cr (3)&\,\,[s_i,s_j] = 0,\,i,j\in I_{d(n-1)}\cr}\eqno (5.89)$$}\vs {\bf Proof.} (1) and (2)
have already been proved. See (4.57) and (5.60). We therefore have only to prove (3). We use the notation of \S 5.6
where we have to prove (5.64). But, as we have observed, this comes down to proving (5.70). Recalling the isomorphism
(5.69) and the definitions of $s_i'$ and $b'$ for $b\in A_{\r}$ in \S 5.6, one has the disjoint union $$W_x = 
\sqcup_{b\in A_{\r}}\,b'\cdot V_{x_o}\eqno (5.90)$$ and $$b'\cdot V_{x_o} = \{y\in V_{x}\mid s_i(y)
=\zeta_i(b^{-1}),\,\forall i\in I_{d(n-1)}\}\eqno (5.91)$$ by (5.87) and (5.88). 

But now, by (5.90), there exists $b\in A_{\r}$ such that $x\in b'\cdot V_{x_o}$. For notational simplicity put $X =
b'\cdot V_{x_o}$. But then $X$
is Lagrangian submanifold of $V_x$ by Propositions 5.19 and 5.20. Obviously $O_x\in {\cal
S}(X)$ (using the notation of Proposition 5.19). Thus $O_x\cap X$ is a Lagrangian submanifold of the symplectic
manifold $O_x$ and $x\in O_x\cap X$. Let $v\in T_x(O_x\cap X)$. But since the $s_i'$ are constant on $X$ by (5.91) 
one has $v\,s_i'= 0$ for all $i\in I_{d(n-1)}$. On the other hand $(\xi_{s_i'})_x\in T_x(O_x)$ by Proposition 1.3.
Thus $\omega_x((\xi_{s_i'})_x,v)= 0$, by (1.17) and (1.19), for all $v$ in the Lagrangian subspace $T_x(O_x\cap X)$ of
$T_x(O_x)$. But by the isotropic maximality of $T_x(O_x\cap X)$, with respect to $\omega_x$, one must have
$(\xi_{s_i'})_x\in T_x(O_x\cap X)$. But then $(\xi_{s_i'})_x\,s_j' = 0$ for all $i,j\in I_{d(n-1)}$ since $s_j'$ is
constant on $O_x\cap X$. This proves (5.70). QED\vs 5.8. The Poisson bracket in $P(n)$ (see \S 1.1) or, as denoted in 
\S 4.4, ${\cal O}(M(n))$, extends in the obvious way to the quotient field, $F(n)$ (see \S 4.5). Similarly, the
Poisson structure in ${\cal O}(M_{\Omega}(n,\ee)$ extends to the function field $F(n,\ee)$ (see \S 4.5) which, we
recall, is a Galois extension of $F(n)$. See Proposition 4.12. 

The Gelfand-Kirillov theorem is the statement that
the quotient division ring of the universal enveloping algebra of $M(n)$ is isomorphic to the quotient division ring of 
a Weyl algebra over a (central) rational function field. A natural commutative analogue is the statement that $F(n)$,
as a Poisson field, is isomorphic to the rational function field of a classical phase space over (a Poisson central)
function field. Using the eigenvalue functions, $r_i$, and the commutative algebraic group $A_{\r}$ we now find that
the statement is explicitly true for the Galois extension $F(n,\ee)$ of $F(n)$.\vs {\bf Theorem 5.24.} {\it For $i\in
I_{d(n-1)}$ one has $s_i^{-1}\in {\cal O}(M_{\Omega}(n,\ee))$ (see (5.58)) so that (see \S 4.2 and (4.47)) 
$r_{(i)}\in {\cal O}(M_{\Omega}(n,\ee))$ where $$r_{(i)} = r_i/s_i\eqno (5.92)$$ Then $F(n,\ee)$ is the rational
function field in $n^2$ variables,
                       $$F(n,\ee) = \Bbb C (r_{(1)},\ldots,r_{(d(n-1))},s_1,\ldots,s_{d(n-1)},r_{d(n-1)
+1},\ldots,r_{d(n-1) + n})\eqno (5.93)$$ Furthermore one has the Poisson commutation relations: $r_{d(n-1) + k}$
Poisson commutes with all every element in $F(n,\ee)$ for $k\in I_n$ and, for $i,j\in I_{d(n-1)}$, $$
\eqalign{(1)&\,\,[r_{(i)},r_{(j)}]= 0\cr (2)&\,\,[r_{(i)},s_j] = \delta_{i\,j}
\cr (3)&\,\,[s_i,s_j] = 0\cr}\eqno (5.94)$$}\vs {\bf Proof.} As a function field one has $$F(n,\ee) = \Bbb C
(r_{1},\ldots,r_{d(n-1)},s_1\ldots,s_{d(n-1)},r_{d(n-1) +1},\ldots,r_{d(n-1) + n})\eqno (5.95)$$ by Theorem 5.14. 
For $i\in
I_{d(n-1)}$ one has $s_i^{-1}\in {\cal O}(M_{\Omega}(n,\ee))$ since $s_i$ vanishes nowhere on $M_{\Omega}(n,\ee)$ by
(5.58). As a function field (5.93) follows immediately from (5.95). The commutation relations (5.94) follow easily
from (5.89). QED\vs Involved in the paper are two groups which operate on
$M_{\Omega}(n,\ee)$ and then contragrediently on ${\cal O}(M_{\Omega}(n,\ee))$. In both cases the latter action
extends to an action on the field $F(n,\ee)$. The first group is the Galois group $\Sigma_n$ (see \S4.5 and
Proposition 4.12). It is immediate from Proposition 4.14 that $\Sigma_n$ preserves Poisson bracket. The second
group is the $d(n-1)$-dimensional complex torus, $A_{\r}$. Since the action of $A_{\r}$ is algebraic (see
Theorem 5.8) the group $A_{\r}$ stabilizes ${\cal O}(M_{\Omega}(n,\ee))$. In fact since the action is algebraic,
as one knows, for any $f\in {\cal O}(M_{\Omega}(n,\ee))$, $$A_{\r}\cdot f\,\,\hbox{spans a finite dimensional
subspace of ${\cal O}(M_{\Omega}(n,\ee))$} \eqno (5.96)$$
In addition, the action of $A_{\r}$ obviously extend to an action on the field
$F(n,\ee)$. Furthermore this action also preserves Poisson bracket (an easy consequence of (5.46) and (5.96)). The
following theorem explicitly determines the action of the two groups. If $\alpha\in \Bbb Z^{d(n-1)}$ and $j\in
I_{d(n-1)}$ let $\alpha_j\in \Bbb Z$ be such that $\alpha = (\alpha_1,\ldots,\alpha_{d(n-1)})$. For $\alpha\in \Bbb
Z$ let $\zeta^{\alpha}$ be the character on the torus $A_{\r}$ defined by putting (see (5.55) and
(5.56)) $$\zeta^{\alpha} = \zeta_1^{\alpha_1}\cdots \zeta_{d(n-1)}^{\alpha_{d(n-1)}}$$ Also let $s^{\alpha}\in
F(n,\ee)$ be defined by putting $$s^{\alpha} = s_1^{\alpha_1}\cdots s_{d(n-1)}^{\alpha_{d(n-1)}}$$ \vskip .5pc
{\bf Theorem 5.25.} {\it Let $m\in I_{n}$ and let $j\in I_{[m]} = I_{d(m)}-I_{d(m-1)}$. Then for any $\sigma\in \Sigma_n$
there exists $k\in I_{[m]}$ such that $$\sigma\cdot r_j = r_k\eqno (5.97)$$ Furthermore if $m\in I_{n-1}$ then $\sigma\cdot
s_j = s_k$ and
$\sigma\cdot r_{(j)} = r_{(k)}$. In addition $\Sigma_n$ normalizes the torus $ A_{\r}$. In fact in the preceding notation $$
\sigma\,A_{\r,j}\,\sigma^{-1} = A_{\r,k}\eqno (5.98)$$ (see \S 5.3 and (5.41). 

Next, for any rational function $f\in \Bbb C(r_1,\ldots,r_{d(n)})$ and any $\alpha\in \Bbb Z^{d(n-1)}$ one has $$b\cdot
(f\,s^{\alpha}) = \zeta^{\alpha}(b)\,f\,s^{\alpha}\eqno (5.99)$$ for any $b\in A_{\r}$.}\vs {\bf Proof.} (5.97) follows
immediately from the definition of the action of $\Sigma_n$. See (4.5), (4.6), (4.10),(4.11), the definition of $\rho_{i\,m}$
in \S 4.2 and (4.47). Since $\sigma$ preserves the Poisson bracket structure in $M_{\Omega}(n,\ee)$ it carries the vector
field $\xi_{r_j}$ to $\xi_{r_k}$ and, hence transforms the corresponding flows (where, by (4.59), we may assume $m\in I_{n-1}$)
so as to yield (5.98). But
$\sigma$ also stabilizes $M_{\Omega}(n,\ee,\b)$ since, by definition, (see \S 5.5) $M_{\Omega}(n,\ee,\b) =
\pi_n^{-1}(\b_{e,\Omega(n)})$, But then $\sigma\cdot s_j = s_k$ by (5.58) and (5.98). It is, of course, then immediate that
$\sigma\cdot r_{(j)} = r_{(k)}$. 
 
The equation (5.99) follows obviously from the definition of $A_{\r}$ (see (5.41)) and (5.59). QED\vskip 1 pc

\centerline{\bf References for Part II}\vskip 1 pc
\item {[B]} A.Borel, {\it Linear Algebraic Groups}, W. A. Benjamin,
Inc, 1969
\item {[GKL1]} A.Gerasimov, S. Kharchev and D. Lebedev, Representation Theory and Quantum Inverse
Scattering Method: The Open Toda Chain and the Hyperbolic Sutherland Model, IMRN 2004 {\bf 17},
2004, 823--854, arXiv:math. QA/0204206
\item {[GKL2]} A.Gerasimov, S. Kharchev and D. Lebedev, On a class of integrable systems connected
with $GL(N,\Bbb R)$, Int. J. Mod. Phys. A, vol 19, Suppl. 2004, 823--854, arXiv: math. QA/0301025
\item {[K2]} B. Kostant, On Whittaker Vectors and Representation
Theory, Inventiones Math. vol 48(1978), 101-184
\item {[KW]} B. Kostant and N. Wallach, Gelfand-Zeitlin theory from
the perspective of classical mechanics I,
http://www.arXiv.org/pdf/math.SG/0408342, to appear in Studies in
Lie Theory, A. Joseph Festschrift, Birkhauser, PM series.
\item {[M1]} D. Mumford, {\it The red book of varieties and
schemes}, Lecture Notes in Math., 1358, Springer, 1995
\item{[M2]} D. Mumford, Geometric Invariant Theory, Ergeb. Math., Band 34(1965), Springer-Verlag
\item {[Sm]} L. Smith, Polynomial Invariants of Finite
Groups, Research Notes in Math. vol 6(1995), A. K. Peters, Ltd
\item {[Sp]} T. Springer, {\it Linear Algebraic Groups}, 2nd
edition, Prog in Math., vol {\bf 9}, Birkhauser Boston, 1998
\item {[ZS]} O. Zariski and P. Samuel,
{\it Commutative Algebra}, Vol 1(1958), Van Nostrand.\vskip 1.5pc
\parindent=30pt
\baselineskip=14pt
$$\vbox{\halign to \hsize{ # \qquad \hfil && # \qquad \hfil\cr
Bertram Kostant && Nolan Wallach \cr
Dept. of Math. && Dept. of Math.\cr
MIT && UCSD \cr
Cambridge, MA 02139 && San Diego, CA 92093\cr
kostant@math.mit.edu && nwallach@ucsd.edu\cr}}
$$
\end

\vskip 1.9pc
\vbox to 60pt{\hbox{Bertram Kostant}
      \hbox{Dept. of Math.}
      \hbox{MIT}
      \hbox{Cambridge, MA 02139}} \noindent E-mail
kostant@math.mit.edu\vskip 1pc
      
\vs\vbox to 60pt{\hbox{Nolan Wallach}
      \hbox{Dept. of Math.}
      \hbox{UCSD}
      \hbox{San Diego, CA 92093}}\noindent
       E-mail nwallach@ucsd.edu

\end

\end